\documentclass{article}
\usepackage[a4paper]{geometry}

\usepackage{amssymb}
\usepackage{amsmath}
\usepackage{amsthm}
\usepackage{mathrsfs}
\usepackage{commath}
\usepackage{mathtools}
\usepackage{dutchcal}

\usepackage{lineno}

\usepackage{stmaryrd}
\usepackage{pifont}
\usepackage{url}
\usepackage{graphicx}
\usepackage[noabbrev]{cleveref}
\usepackage{csquotes}
\usepackage{adjustbox}
\usepackage{textcomp}

\usepackage{multirow}

\usepackage[round,sort,authoryear,semicolon]{natbib}

\usepackage{algorithm}
\usepackage{algpseudocode}

\algdef{SE}[DOWHILE]{Do}{doWhile}{\algorithmicdo}[1]{\algorithmicwhile\ #1}%
\algnewcommand\algorithmicforeach{\textbf{for each}}
\algdef{S}[FOR]{ForEach}[1]{\algorithmicforeach\ #1\ \algorithmicdo}

\usepackage{authblk}

\usepackage{siunitx}

\usepackage{pgfplots}
\usepackage{booktabs}
\usepgfplotslibrary{polar}
\usetikzlibrary{calc}

\newcommand{\avg}[1]{\ensuremath{\{\!\!\{{#1}\}\!\!\}}}
\newcommand{\jump}[1]{\llbracket{#1}\rrbracket}
\newcommand{\tjump}[1]{\llbracket{#1}\rrbracket_{\times}}
\newcommand{\ojump}[1]{\llbracket{#1}\rrbracket_{\otimes}}
\newcommand{\dgnorm}[1]{\|\hspace{-0.2ex}| {#1} \|\hspace{-0.2ex}|_\text{IP}}
\newcommand{\snorm}[1]{\left\lvert {#1} \right\rvert}

\newcommand{\avga}[1]{\ensuremath{\avg{#1}_{\ell}}}
\newcommand{\avgmp}[1]{\ensuremath{\avg{#1}_{\ell^\mp}}}

\renewcommand\vec{\boldsymbol}

\DeclareSIUnit\year{yr}

\newtheorem{theorem}{Theorem}[section]

\newtheorem{lemma}[theorem]{Lemma}
\newtheorem{remark}[theorem]{Remark}

\usepackage[acronym]{glossaries-prefix}
\makenoidxglossaries

\newacronym{FE}{FE}{finite element}
\newacronym{EDG}{EDG}{embedded discontinuous Galerkin}
\newacronym{DG}{DG}{discontinuous Galerkin}
\newacronym{HDG}{HDG}{hybrid discontinuous Galerkin}
\newacronym{TH}{TH}{Taylor--Hood}
\newacronym{PDE}{PDE}{partial differential equation}
\newacronym{ODE}{ODE}{ordinary differential equation}
\newacronym{RK}{RK}{Runge--Kutta}
\newacronym{CFL}{CFL}{Courant--Friedrichs--Lewy}
\newacronym{ALA}{ALA}{anelastic liquid approximation}
\newacronym{IP}{IP}{interior penalty}
\newacronym{UFL}{UFL}{Unified Form Language}
\newacronym{FFC}{FFC}{FEniCS Form Compiler}
\newacronym[prefixfirst={the\ }]{PETSc}{PETSc}{Portable, Extensible Toolkit for Scientific Computation}
\newacronym{CAD}{CAD}{computer--aided design}
\newacronym{DoF}{DoF}{degree of freedom}
\newacronym[prefixfirst={a\ }, prefix={an\ }]{FGMRES}{FGMRES}{flexible generalized minimal residual}
\newacronym[prefixfirst={a\ }, prefix={a\ }]{GMRES}{GMRES}{generalized minimal residual}
\newacronym{ILU}{iLU}{incomplete LU}
\newacronym{SIPG}{SIPG}{symmetric interior penalty Galerkin}
\newacronym{NIPG}{NIPG}{nonsymmetric interior penalty Galerkin}
\newacronym{IPG}{IPG}{interior penalty Galerkin}
\newacronym{RIPG}{RIPG}{robust interior penalty Galerkin}

\usepackage{enumitem}
\crefname{enumi}{caveat}{caveats}
\Crefname{enumi}{Caveat}{Caveats}

\let\oldequation\equation
\let\oldendequation\endequation
\renewenvironment{equation}
  {\linenomathNonumbers\oldequation}
  {\oldendequation\endlinenomath}

\let\oldalign\align
\let\oldendalign\endalign
\renewenvironment{align}
  {\linenomathNonumbers\oldalign}
  {\oldendalign\endlinenomath}

\let\oldalignat\alignat
\let\oldendalignat\endalignat
\renewenvironment{alignat}
  {\linenomathNonumbers\oldalignat}
  {\oldendalignat\endlinenomath}

\begin{document}

\title{A divergence free $C^0$-RIPG stream function formulation of the
incompressible Stokes system with variable viscosity}

\author[1,*]{Nathan Sime}
\author[2]{Paul Houston}
\author[1]{Cian R. Wilson}
\author[1]{Peter E. van Keken}
\affil[1]{Earth and Planets Laboratory, Carnegie Institution for Science, Washington, D.C., U.S.A.}
\affil[2]{School of Mathematical Sciences, University of Nottingham, Nottingham, U.K.}
\affil[*]{Corresponding author: Nathan Sime, nsime@carnegiescience.edu}

\maketitle

\begin{abstract}
Pointwise divergence free velocity field approximations of the Stokes system
are gaining popularity due to their necessity in precise modelling of physical
flow phenomena. Several methods have been designed to satisfy this requirement;
however, these typically come at a greater cost when compared with standard
conforming methods, for example, 
because of the
complex implementation and development of specialized
finite element bases. Motivated by the desire to mitigate these issues for 2D
simulations,
we present a $C^0$-\gls{IPG} discretization 
of the Stokes system
in the stream function formulation.
In order to preserve a spatially varying
viscosity this approach does not yield the standard and well known biharmonic
problem. We further employ the so-called \gls{RIPG} method; stability and convergence
analysis
of the proposed scheme is undertaken. The former, which involves deriving a bound on the 
interior penalty parameter is particularly
useful to address the $\mathcal{O}(h^{-4})$ growth in the condition number of
the discretized operator. Numerical experiments confirming the optimal convergence of the 
proposed method are undertaken. Comparisons with thermally driven buoyancy mantle convection model benchmarks
are presented.
\end{abstract}

\section{Introduction}

\subsection{Motivation}

Large areas of computational fluid dynamics rely on the accurate solution of the Stokes system which
requires correct satisfaction of the mass conservation equation. 
We are particularly interested in geophysical flow calculations that model slow convection in the 
Earth's mantle to help understand the Earth's thermal and chemical evolution 
\citep{Schubert2001,Bercovici2015,Ricard2015}.
Specific applications include: the modeling of the recycling and subsequent mixing of oceanic crust 
\citep{christensen1994,brandenburg2008}
potentially within preexisting heterogeneity 
\citep{Gulcher2021},
the thermal and chemical evolution of subduction zones 
\citep{WadaKing2015,Gerya2021}; 
and the formation of hotspots and large igneous provinces by mantle plumes.

\begin{figure}[t!]
\centering
\includegraphics[width=0.6\linewidth]{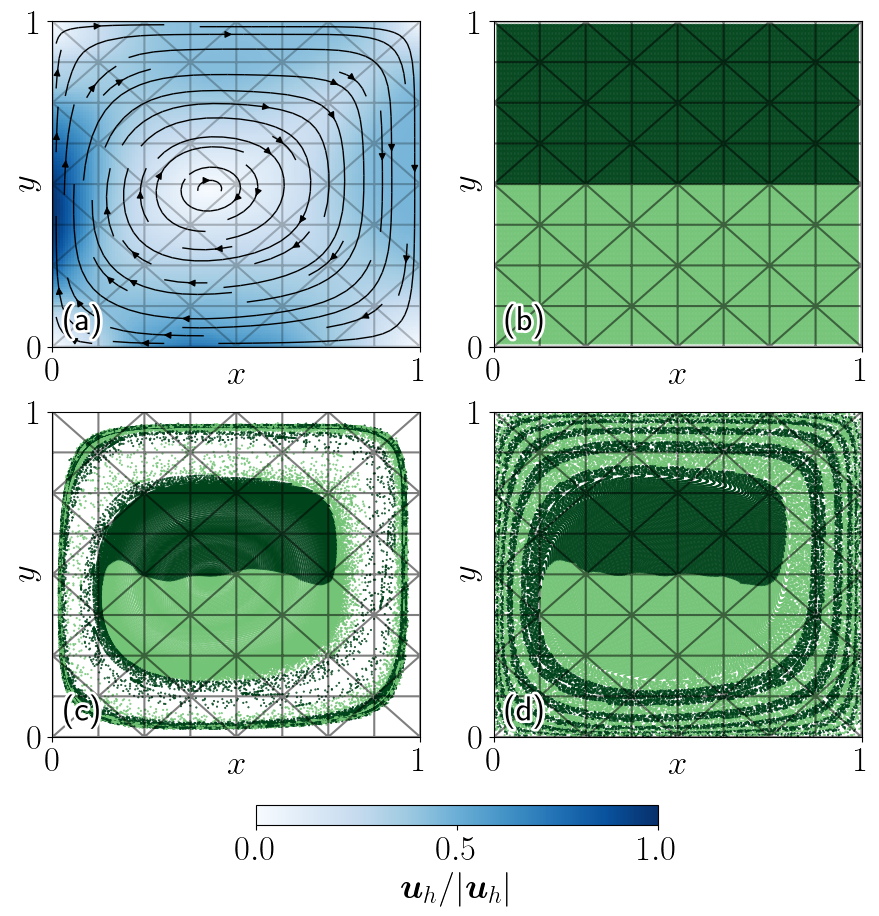}
\caption{Advecting particles through velocity fields which are not pointwise
divergence free may lead to spurious results. Here, we show an example
advecting tracers using a third order accurate Runge--Kutta scheme in a mantle
convection model exhibited as case~T4 in our numerical experiments later in
\Cref{sec:tosi_problem} (comprising a numerical benchmark in
\citet{tosi2015}). (a)~The velocity field approximation computed using the
$C^0$-RIPG method developed in this work; shown are the computational mesh,
velocity streamlines and normalized speed field. (b)~The initial configuration
of $256^2$ equidistant particles colored solely as a visual aid. (c)~The
advected particles after approximately 10~overturns in a \emph{non}-pointwise
divergence free velocity field computed using the \gls{TH} method. (d)~The
advected particles after approximately 10~overturns in the pointwise
divergence free velocity field computed using the $C^0$-RIPG method developed
in this work. The average number of particles per cell is $512$ with standard
deviation of $36.67$, $240.06$ and $35.00$ in (b), (c) and (d), respectively.
The velocity field is discretized using polynomials of degree 2 in both the
\gls{TH} and $C^0$-RIPG cases shown.}
\label{fig:tracer_intro}
\end{figure}

Exact, or at least approximately local, mass conservation is particularly important for particle methods 
\citep{christensen1994,vanKeken1997,tackley2003}.
It has been shown, for example, that
significant artefacts can occur when using the popular \gls{TH} finite element method
for the numerical approximation of
the Stokes system with particles. The \gls{TH} element pair for velocity-pressure does not satisfy mass conservation
locally. Such artefacts include the formation of holes and concretions in the particle distribution
or the artificial settling of particles in gravity driven flows. Recent geodynamical examples demonstrating this are for 
purely compositionally-driven flow 
\citep{samuel2018,sime2020b}, 
thermochemical convection 
\citep{jones2021,Wang2015,pusok2017}, and when mathematical fields represented by particles are advected 
\citep{maljaars2019,maljaars2020,sime2020b,sime2022}. We demonstrate these artifacts in \Cref{fig:tracer_intro}
and further refer to \citet{Jenny2001,McDermott2008}.

In this paper we will explore a new approach to guarantee exact mass conservation. Here, we specifically 
consider the incompressible Stokes system in 2D in a simply connected domain;
this covers many of the typical geometries that are frequently used in geodynamical applications 
\citep{brandenburg2008,Hernlund2008,jones2021,LiMcNamara2022}.
The extension to simple compressible flow modeling for mantle convection 
\citep{jarvis_mckenzie_1980,Tackley2008,bossmann2013}
is possible in the same fashion as in the transition made from \cite{sime2020b} to \cite{sime2022}.

\subsection{Context}

Several numerical schemes have been developed to exactly satisfy mass
conservation of the velocity approximation, denoted here by $\vec{u}_h$, in a pointwise sense,
by which we mean that $\nabla \cdot \vec{u}_h(\vec{x}) = 0$ for all $\vec{x}$
in the computational
domain~\citep[e.g.,][]{brezzi1985,cockburn2011analysis,volker2017,raviart1977,rhebergen2018hybridizable,rhebergen2020,scott1985}.

In our recent work, which demonstrated the need for exact mass conservation 
~\citep{jones2021,sime2020b,sime2022}, 
we employed an
embedded discontinuous Galerkin--\gls{HDG} method which yields a solenoidal velocity
approximation (as developed in \citet{rhebergen2020} and implemented in
\citet{maljaars2020}). However, even with the utilization of static condensation, which
factorizes the element-local problems in favor of the global system defined on
the facets, the system assembly and computation of the solution is expensive. Seeking a more
computationally efficient numerical scheme in the 2D setting is the key motivation 
of this article.

The general approach which the above mentioned numerical schemes exploit to achieve exact
pointwise satisfaction of mass conservation is by ensuring that the divergence
of the velocity approximation lies in the space in which we seek the pressure
approximation, denoted by $Q^h$, i.e., that $\nabla \cdot \vec{u}_h \in Q^h$. 
These schemes typically require the definition of specialized \gls{FE} basis functions which may be difficult
to implement. Libraries such as \texttt{basix}~\citep{scroggs2022} simplify and
automate the creation of these \gls{FE} bases; however, the variational
formulation, assembly and computational solution of these systems may remain
difficult. Also, as previously mentioned above, assembling \gls{HDG} systems via static
condensation requires careful management to preserve scalable and efficient
local assembly into a global system.

A popular alternative is to consider the biharmonic formulation of the Stokes system which
provides an exactly divergence free approximation of the velocity. In
this setting, the velocity is cast as the curl of an unknown potential field, i.e.,
$\vec{u}_h = \nabla \times \phi_h$, where $\phi_h$ is the \gls{FE} approximation of the 
stream function. Hence, the numerical approximation of the velocity vector is trivially solenoidal.
The standard \gls{FE} discretization of the fourth-order biharmonic problem
requires a $C^1$ conforming basis for the approximation of $\phi_h$
\citep[e.g.,][]{morley1968triangular,argyris1968tuba} or, for example, a $C^k$, $k \geq 1$,
continuous divergence conforming B-spline basis
\citep[e.g.,][]{christensen1994,evans2013isogeometric,kopitzke1979,vanKeken1993}.

In this work, our objective is to develop a \gls{FE} formulation of the Stokes
system by taking the same approach as the $C^0$-interior penalty Galerkin (\gls{IPG}) formulation of the
biharmonic problem, whilst preserving the symmetric gradient operator
permitting variable viscosity models. \gls{IPG} methods in the context of
\gls{DG} \gls{FE} discretizations for second order \glspl{PDE} are well
understood~\citep[see, for example,][]{unifieddg}. In the context of fourth-order
problems we seek to extend the technique presented in, for example,
\citet{engel2002}; see also the recent complimentary paper by \citet{dong_mascotti_2023}. In this setting, $\phi_h$ is approximated with a $C^0$~\gls{FE} basis and continuity of its first derivative is enforced
weakly. This formulation will yield many new terms on the facets of the
\gls{FE} mesh which may be verbose and arduous to implement. However, we
exploit modern computational symbolic specification \citep{ufl} and automatic
generation of \gls{FE} formulations \citep{FFC} which vastly reduces implementation
complexity.

Given this objective, we highlight what we see as clear advantages of the
$C^0$-\gls{IPG} formulation and therefore our motivation:
\begin{enumerate}
    \item Solenoidal velocity approximation.
    \item Exploiting strong imposition of $\phi_h = 0$ on the computational
    boundary yields pointwise satisfaction of $\vec{u}_h \cdot \vec{n} = 0$,
    where $\vec{n}$ denotes the unit outward pointing normal vector on the boundary of the computational domain.
    \item In the context of the $2D$ system, reduction of the Stokes system's
    unknown velocity and pressure variables to a single unknown scalar
    potential.
    \item Assembly using a standard $C^0$~\gls{FE} basis with no special
    requirements.
\end{enumerate}
However, we must highlight that the $C^0$-\gls{IPG} method does come with
caveats:
\begin{enumerate}
    \item The condition number of the underlying matrix stemming from the
    discretization of the fourth-order \gls{PDE} exhibits
    growth at a rate of approximately $h^{-4}$ as the mesh is uniformly refined; here, $h$ is a measure of the mesh
    element size. \label{pt:caveat_h4}
    \item A dimensionless interior penalty parameter must be chosen to ensure
    stability of the scheme. \label{pt:caveat_penalty}
    \item The matrix associated with the $C^0$-\gls{IPG} discretization is more
    dense than a standard conforming method.
\end{enumerate}
\Cref{pt:caveat_penalty} further impacts \cref{pt:caveat_h4} in the sense that the
interior penalty parameter should be carefully selected not to be so large as to
further adversely affect the condition number of the underlying matrix. For this reason,
we employ the robust \gls{IPG} (\gls{RIPG}) formulation introduced in \citet{dong2022}. 
Here, the key modification to the underlying \gls{IPG} scheme is the exploitation of
a weighted average operator; this then allows for a lower bound on the interior
penalty parameter to be determined in a simple manner.

\subsection{Structure of the remainder of the paper}

This article is structured as follows: in \Cref{sec:streamfunction} we define the
Stokes system and the associated stream function formulation. In \Cref{sec:c0ipg}
we define the $C^0$-\gls{RIPG} formulation, based on employing a suitable weighted average operator
and interior penalty parameter. The stability and error analysis of the 
proposed scheme is studied in \Cref{sec:stability_and_error_analysis} where $hp$--optimal
error bounds are established with respect to a given DG energy norm; this is in accordance
with the analogous results derived in the case of a constant viscosity coefficient \citep{dong_mascotti_2023}.
The numerical performance of the
$C^0$-\gls{RIPG} scheme is considered in \Cref{sec:experiments}. 
Finally, in \Cref{sec:conclusion} we summarize the work presented in this article and 
discuss potential future developments.

\section{Stream function formulation}
\label{sec:streamfunction}

Let $\Omega \subset \mathbb{R}^2$ denote a simply connected domain with boundary
$\partial \Omega$ and outward pointing unit normal vector $\vec{n}$. The boundary is
subdivided into Dirichlet and Neumann components $\partial\Omega_D$ and
$\partial\Omega_N$, respectively, which do not overlap, i.e.,
$\partial\Omega_D \cup \partial\Omega_N = \partial\Omega$ and
$\partial\Omega_D \cap \partial\Omega_N = \emptyset$.
In $\Omega$ we seek the velocity $\vec{u} : \Omega \mapsto
\mathbb{R}^2$ and pressure $P : \Omega \mapsto
\mathbb{R}$ which satisfy the Stokes system
\begin{align}
-\nabla \cdot (2 \mu \varepsilon(\vec{u})) + \nabla P &= \vec{f}, \label{eq:momentum} \\
\nabla \cdot \vec{u} &= 0, \label{eq:mass}
\end{align}
subject to the boundary conditions
\begin{alignat}{2}
\vec{u} &= \vec{0} &&\text{ on } \partial\Omega_D, \label{eq:bcs-1} \\
2 \mu \varepsilon(\vec{u}) \cdot \vec{n} - P \vec{n} &= \vec{g}_N &&\text{ on } \partial\Omega_N. \label{eq:bcs-2}
\end{alignat}
Here, $\varepsilon(\vec{u}) = \frac{1}{2}(\nabla \vec{u} + \nabla
\vec{u}^\top)$ is the symmetric gradient, $\mu(\vec{x}) : \Omega \mapsto \mathbb{R}^+$
is the viscosity, and $\vec{f}(\vec{x}) : \Omega \mapsto \mathbb{R}^2$ is a given forcing function.

Defining the space $\vec{H}^1_{\vec{0}}(\Omega) = \left\{ \vec{v} \in [H^1(\Omega)]^2: \vec{v} = \vec{0} \text{ on } \partial\Omega_D \right\}$,
%
%
the weak formulation of \cref{eq:momentum,eq:mass} reads: find $(\vec{u}, P) \in \vec{H}^1_{\vec{0}}(\Omega) \times L_2(\Omega)$ such that
\begin{align}
B(\vec{u}, \vec{v}) - (\nabla \cdot \vec{v}, P) &= l(\vec{v}), \\
(\nabla \cdot \vec{u}, q) &= 0
\end{align}
for all $(\vec{v}, q) \in \vec{H}^1_{\vec{0}}(\Omega) \times L_2(\Omega)$, where
\begin{align}
B(\vec{u}, \vec{v}) &= (2 \mu \varepsilon(\vec{u}), \varepsilon(\vec{v})), \\
l(\vec{v}) &= (\vec{f}, \vec{v}) + \langle \vec{g}_N, \vec{v} \rangle_{\partial\Omega_N}.
\end{align}
Here, $(\cdot,\cdot)$ is the standard $L_2(\Omega)$ inner product and $\langle \cdot,\cdot \rangle_{\omega}$ denotes the $L_2(\omega)$ inner product on a subset $\omega$ of the boundary $\partial\Omega$.

Defining the divergence free space $\vec{Z}_{\vec{0}} =
\{\vec{v} \in \vec{H}^1_{\vec{0}}(\Omega) : \nabla \cdot
\vec{v} = 0 \}$ allows us to state the divergence free weak formulation:
find $\vec{u} \in \vec{Z}_{\vec{0}}$ such that
\begin{equation}
B(\vec{u}, \vec{v}) = l(\vec{v})
\quad\quad \forall \vec{v} \in \vec{Z}_{\vec{0}}.
\end{equation}
Given that $\Omega \subset \mathbb{R}^2$ is simply connected, then for every $\vec{v} \in \vec{Z}_{\vec{0}}$, there exists a unique $\psi \in H^2(\Omega) \setminus \mathbb{R}$ such that
$\vec{v} = \nabla \times \psi$, where we define the curl vector acting on a scalar, $\psi : \Omega \mapsto \mathbb{R}$
and a vector, $\vec{z} = (\vec{z}_x, \vec{z}_y)^\top : \Omega \mapsto \mathbb{R}^2$, by
\begin{equation}
\nabla \times \psi =
\left(\frac{\partial \psi}{\partial y}, - \frac{\partial \psi}{\partial x}\right)^\top
\quad \text{and} \quad
\nabla \times \vec{z} =
\frac{\partial \vec{z}_y}{\partial x}
- \frac{\partial \vec{z}_x}{\partial y},
\end{equation}
respectively, cf.~\citet{gir_rav_ns_86}.
Setting $\Psi_{\vec{0}} = \{ \psi \in H^2(\Omega) \setminus
\mathbb{R}: \nabla \times \psi = \vec{0} \text{ on } \partial\Omega_D \}$
we write the stream function formulation: find $\phi \in \Psi_{\vec{0}}$
such that
\begin{equation}
B(\nabla \times \phi, \nabla \times \psi) = l(\nabla \times \psi)
\quad \quad \forall \psi \in \Psi_{\vec{0}}.
\label{eq:stream_function_formulation}
\end{equation}


\section{$C^0$-\gls{IPG} formulation}
\label{sec:c0ipg}

Let $\mathcal{T}^h$ be the subdivision of $\Omega$ into a mesh composed of a
tessellation of nonoverlapping triangular elements $\kappa$ such that $\mathcal{T}^h = \{
\kappa \}$ and $\overline{\Omega} = \cup_{\kappa\in\mathcal{T}^h}
\overline{\kappa}$. Each element has boundary $\partial\kappa$ with outward
pointing unit normal vector $\vec{n}_\kappa$. The interior facets of the mesh are
defined by $\Gamma_I = \cup_{\kappa\in\mathcal{T}^h} \partial \kappa \setminus
\partial\Omega$, the exterior Dirichlet facets by $\Gamma_E =
\cup_{\kappa\in\mathcal{T}^h} \partial \kappa \cap \partial \Omega_D$ and the mesh
skeleton by $\Gamma = \Gamma_I \cup \Gamma_E$. Given two neighboring elements
$\kappa^+$ and $\kappa^-$ and the smooth functions $\phi^\pm : \kappa^\pm
\mapsto \mathbb{R}$ and $\vec{v}^\pm : \kappa^\pm \mapsto
\mathbb{R}^2$, we define on their common face $F = \partial\kappa^+ \cap
\partial\kappa^-$ the following operators
\begin{align}
\avg{\phi} &= w^+ \phi^+ + w^- \phi^-, \\
\ojump{\vec{v}} &= \vec{v}^+ \otimes \vec{n}^+ + \vec{v}^- \otimes \vec{n}^-,
\end{align}
which we refer to as the weighted average and tensor jump, respectively; for notational
simplicity, here we write $\vec{n}^\pm=\vec{n}_{\kappa^\pm}$, respectively. On a boundary face
$F = \partial\kappa \cap \partial\Omega_D$, we set $\avg{\phi} = \phi$ and
$\ojump{\vec{v}} = \vec{v} \otimes \vec{n}$, where $\vec{n}$ denotes the outward unit
normal vector on the boundary $\partial\Omega_D$.
Here, the
weights $w^\pm$ are to be determined; however, they are subject to the
constraint that $w^+ + w^- = 1$. Furthermore, we define the broken divergence,
gradient, curl and symmetric gradient operators on each $\kappa
\in \mathcal{T}^h$
\begin{equation}
(\nabla_h \cdot (\cdot))|_{\kappa} = \nabla \cdot (\cdot|_\kappa), \quad
(\nabla_h (\cdot))|_{\kappa} = \nabla (\cdot|_\kappa), \quad
(\nabla_h \times (\cdot))|_{\kappa} = \nabla \times (\cdot|_\kappa), \quad
\varepsilon_h(\cdot) = \frac{1}{2}\left(\nabla_h (\cdot) + \nabla_h(\cdot)^\top \right),
\end{equation}
respectively. 


We define the scalar \gls{FE} function space
\begin{equation}
V^{h,p}_0 = \{ v \in H^1_0(\Omega) : \left. v \right|_{\kappa} \in \mathcal{P}_p(\kappa) \:\: \forall \kappa \in \mathcal{T}^h\},
\end{equation}
where $H^1_{0}(\Omega) = \left\{ v \in H^1(\Omega): v = 0 \text{ on }
\partial\Omega_D \right\}$, $\mathcal{P}_p(\kappa)$ denotes the polynomials of
total degree $p$ on $\kappa$, $\kappa \in \mathcal{T}^h$, and $p \geq 2$. Note
that $V^{h,p}_0 \not\subseteq \Psi_{\vec{0}}$; hence, we will weakly impose the
continuity of the derivatives of the underlying basis by employing the
following $C^0$-\gls{IPG} scheme: find $\phi_h \in V^{h,p}_0$ such that
\begin{equation}
B_\text{IP}(\vec{u}_h, \vec{v}_h) = l_\text{IP}(\vec{v}_h) \quad \forall \psi_h \in V^{h,p}_0,
\label{eq:bilinear_fem}
\end{equation}
where $\vec{u}_h = \nabla_h \times \phi_h$, $\vec{v}_h = \nabla_h \times \psi_h$,
%
%
\begin{align}
B_\text{IP}(\vec{u}, \vec{v}) =
B_h(\vec{u}, \vec{v})
&- \langle \ojump{\vec{u}}, \avg{2 \mu \varepsilon_h(\vec{v})} \rangle_{\Gamma} 
- \langle \avg{2 \mu \varepsilon_h(\vec{u})}, \ojump{\vec{v}} \rangle_{\Gamma}
%
+ \langle \beta \ojump{\vec{u}}, \ojump{\vec{v}} \rangle_{\Gamma},
%
\label{eq:bilinear_B}
\end{align}
$B_h(\vec{u}, \vec{v}) = (2 \mu \varepsilon_h(\vec{u}), \varepsilon_h(\vec{v}))$ and
%
%
$l_\text{IP}(\cdot) = l(\cdot)$.
Here, $\beta$ is the interior penalty parameter which will be defined in the following
section. A derivation of this $C^0$-\gls{IPG} formulation is outlined in
\Cref{sec:appdx_derivation}.

\begin{remark}
In the case when inhomogeneous Dirichlet boundary conditions are employed, i.e., $\vec{u} = \vec{u}_D$ on $\partial\Omega_D$, we seek $\phi_h \in V^{h,p} = \{ v \in H^1(\Omega) : \left.v\right|_\kappa \in \mathcal{P}_p(\kappa) \:\: \forall \kappa \in \mathcal{T}^h \}$ such that
\begin{equation}
B_\text{IP}(\vec{u}_h, \vec{v}_h) = l_\text{IP}(\vec{v}_h) \quad \forall \psi_h \in V^{h,p},
\end{equation}
where $l_\text{IP}(\cdot)$ is replaced by the alternative linear functional
\begin{equation}
l_\text{IP}(\vec{v}) = l(\vec{v})
+ \langle \vec{u}_D \otimes \vec{n}, \beta \vec{v} \otimes \vec{n} - 2 \mu \varepsilon_h(\vec{v}) \rangle_{\partial\Omega_D}.
\label{eq:bilinear_l}
\end{equation}
\end{remark}

\begin{remark}
In the case of free slip and zero penetration conditions we
divide $\partial\Omega_D$ into free slip and zero penetration components
$\partial\Omega_\text{FS}$ and $\partial\Omega_\text{ZP}$, respectively, such that
$\partial\Omega_D = \partial\Omega_\text{FS} \cup \partial\Omega_\text{ZP}$ and
$\partial\Omega_\text{FS} \cap \partial\Omega_\text{ZP} = \emptyset$. On these
components we impose
\begin{align}
\vec{u} \cdot \vec{n} &= 0 &&\text{ on } \partial\Omega_\text{FS} \cup \partial\Omega_\text{ZP}, \label{eq:bcs-3} \\
\left(2 \mu \varepsilon_h(\vec{u}) \cdot \vec{n} - P \vec{n}\right) \cdot \vec{\tau} &= 0 &&\text{ on } \partial\Omega_\text{FS}. \label{eq:bcs-4} \\
\vec{u} \cdot \vec{\tau} &= \vec{u}_D \cdot \vec{\tau} &&\text{ on } \partial\Omega_\text{ZP}, \label{eq:bcs-5}
\end{align}
where $\vec{\tau}$ is the unit vector which lies tangential to
$\partial\Omega$. We may exploit the strong imposition of
$\left.\phi_h\right|_{\partial\Omega_D} = 0$ by seeking $\phi_h \in V^{h,p}_0$
such that
\begin{equation}
B_\text{IP}(\vec{u}_h, \vec{v}_h) = l_\text{IP}(\vec{v}_h) \quad \forall \psi_h \in V^{h,p}_0,
\label{eq:bilinear_fem_freeslip}
\end{equation}
where $l_\text{IP}(\cdot)$ is the alternative linear functional
\begin{equation}
l_\text{IP}(\vec{v}) = l(\vec{v})
+ \langle \vec{u}_D \otimes \vec{n}, \beta \vec{v} \otimes \vec{n} - 2 \mu \varepsilon_h(\vec{v}) \rangle_{\partial\Omega_\text{ZP}}.
\label{eq:bilinear_l_fs_zp}
\end{equation}
This scheme yields pointwise satisfaction of the boundary flux such that
\begin{equation}
\Vert \nabla_h \phi_h \cdot \vec{\tau} \Vert_{L_2(\partial\Omega_D)}
= \Vert \nabla_h \times \phi_h \cdot \vec{n} \Vert_{L_2(\partial\Omega_D)} 
= \Vert \vec{u}_h \cdot \vec{n} \Vert_{L_2(\partial\Omega_D)} 
= 0.
\end{equation}
\end{remark}


\subsection{$C^0$-\gls{RIPG} formulation}
\label{sec:c0ripg}

The classical \gls{SIPG} formulation sets $w^\pm = \frac{1}{2}$. However we seek
to exploit the \gls{RIPG} formulation introduced in \citet{dong2022} which is based
on a particular choice of the weights $w^\pm$ and the interior penalty parameter $\beta$. 
The key advantage of the \gls{RIPG} formulation is
that it remains parameterless in the sense that there is a known bound on $\beta$
which ensures stability of the numerical scheme given appropriate choices of
$w^\pm$. We first state the values we employ in the context of our 2D numerical
experiments when $\mathcal{T}^h$ is composed of triangles. In the
remainder of this section we rationalize this choice based on studying the stability
of the \gls{RIPG} scheme, along with discussing the
advantages when using this method when compared with a standard
\gls{SIPG} numerical scheme.

The \gls{RIPG} scheme average operator weights and interior penalty parameter are 
defined, respectively, by
\begin{equation}
\label{eq_weights}
w^\pm = \frac{\zeta^\pm}{\zeta^+ + \zeta^-},
\quad
\left. \beta \right|_{F} =
\begin{cases}
(\zeta^+ + \zeta^-)^{-2} & F \subset \Gamma_I, \\
(\zeta^+)^{-2} & F \subset \Gamma_E, 
\end{cases}
\end{equation}
where
\begin{equation} \label{eq_zeta}
\left.\zeta^\pm\right|_{F} = \left(
  \delta \sqrt{\frac{3 p(p-1)}{2}  \frac{\snorm{F}}{\snorm{\kappa^\pm}}}
  \norm{2 \mu \vec{n}^\pm}_{L_\infty(F)} 
  \norm{(2 \mu)^{-\frac{1}{2}}}_{L_\infty(\kappa^\pm)}
\right)^{-1},
\end{equation}
where for a set $\omega \subset {\mathbb R}^n$, $n\geq 1$, we write $|\omega|$
to denote the $n$--dimensional Hausdorff measure of $\omega$ and $\delta >
\sqrt{2}$ is a constant, cf. below.

\section{Stability and error analysis}
\label{sec:stability_and_error_analysis}

The aim of this section is to study the stability of the $C^0$-\gls{RIPG} formulation
\eqref{eq:bilinear_fem} and establish an optimal $hp$--error bound. For the
purposes of the proceeding error analysis we introduce a suitable extension of the
bilinear form $B_\text{IP}(\cdot, \cdot)$. To this end, we write 
$\Pi_{L_2}:[L_2(\Omega)]^{2\times 2} \mapsto [V^{h,p-2}]^{2\times 2}$ to denote the
orthogonal $L_2$-projection operator. With this notation for
$\phi, ~\psi \in V:= H^2(\Omega)+V^{h,p}_0$ we write
\begin{align}
\tilde{B}_\text{IP}(\vec{u}, \vec{v}) =
B_h(\vec{u}, \vec{v})
&- \langle \ojump{\vec{u}}, \avg{2 \mu \Pi_{L_2}(\varepsilon_h(\vec{v}))} \rangle_{\Gamma} 
- \langle \avg{2 \mu \Pi_{L_2}(\varepsilon_h(\vec{u}))}, \ojump{\vec{v}} \rangle_{\Gamma}
%
+ \langle \beta \ojump{\vec{u}}, \ojump{\vec{v}} \rangle_{\Gamma},
%
\label{eq:extended_bilinear_B}
\end{align}
where $\vec{u} = \nabla_h \times \phi$ and $\vec{v} = \nabla_h \times \psi$. Furthermore, we introduce
the following DG norm
\begin{equation}
\dgnorm{\vec{v}}^2 = \norm{\sqrt{2 \mu} \varepsilon_h(\vec{v})}^2_{L_2(\Omega)}
+ \norm{\sqrt{\beta} \ojump{\vec{v}}}^2_{L_2(\Gamma)}.
\label{eq:ipnorm}
\end{equation}

\subsection{Coercivity and continuity} \label{sec:coercivity}

The aim of this section is to study the stability of the $C^0$-\gls{RIPG} scheme
\eqref{eq:bilinear_fem}. To this end, we first recall the following inverse equality
from \citet{warburton2003}.

\begin{lemma}
\label{lem:inverse_inequality}
Given $\kappa\in \mathcal{T}^h$ is a triangular element in 2D ($d=2$ below), let $F\subset \partial\kappa$ denote one
of its faces. Then for $v\in\mathcal{P}_p(\kappa)$ the following inverse inequality holds
\begin{align}
\Vert v \Vert_{L_2(F)} &\leq C_\text{inv}(\kappa, F, p) \Vert v \Vert_{L_2(\kappa)},
\end{align}
where
\begin{align}
C_\text{inv}(\kappa, F, p) &= \sqrt{\frac{(p + 1)(p + d)}{d} \frac{|F|}{|\kappa|}}.
\end{align}
\end{lemma}

Equipped with Lemma~\ref{lem:inverse_inequality} we now state the main result of this section.
\begin{lemma}
The bilinear form $\tilde{B}_\text{IP}(\cdot, \cdot)$ is coercive for any $\delta>\sqrt{2}$ and continuous over $V\times V$; in particular
we have that
\begin{align}
\tilde{B}_\text{IP}(\vec{v}, \vec{v}) &\geq C_{\text{coer}}  \dgnorm{\vec{v}}^2 & \forall \phi \in V, \nonumber \\
\tilde{B}_\text{IP}(\vec{v}, \vec{w}) &\leq C_{\text{cont}}  \dgnorm{\vec{v}} \dgnorm{\vec{w}} & \forall \phi,\psi \in V, \nonumber
\end{align}
where $\vec{v} = \nabla_h\times \phi$, $\vec{w} = \nabla_h\times \psi$, $C_{\text{coer}} = 1-\sqrt{2}/\delta$ and $C_{\text{cont}} = 2(1+2/\delta^2)$.
\end{lemma}
\begin{proof}
We first consider the coercivity of the bilinear form $\tilde{B}_\text{IP}(\cdot, \cdot)$. For $\phi \in V$, 
$\vec{v} = \nabla_h\times \phi$, we note that
\begin{equation}
\tilde{B}_\text{IP}(\vec{v}, \vec{v}) = 
\norm{\sqrt{2 \mu} \varepsilon_h(\vec{v})}^2_{L_2(\Omega)}
+ \norm{\sqrt{\beta} \ojump{\vec{v}}}^2_{L_2(\Gamma)}
- 2 \int_\Gamma \avg{2 \mu \Pi_{L_2}(\varepsilon_h(\vec{v}))} \cdot \ojump{\vec{v}} \dif s.
\label{eq:bilin_anal}
\end{equation}
In order to determine the lower bound on $\delta$, we must express the last term in
\cref{eq:bilin_anal} in terms of the first two. Exploiting the inverse inequality stated in Lemma~\ref{lem:inverse_inequality}, together with the
stability of the $L_2$-projection operator $\Pi_{L_2}$, for interior faces, we deduce that
\begin{align}
&\left| 2 \int_{\Gamma_I} \avg{2 \mu \Pi_{L_2}(\varepsilon_h(\vec{v}))} \cdot \ojump{\vec{v}} \dif s \right| \nonumber \\
&\quad\quad\leq
2 \sum_{F \in \Gamma_I} \int_F  \left( \sum_{*\in\{+,-\}} w^* \norm{2 \mu \vec{n}^*}_{L_\infty(F)} \snorm{\Pi_{L_2}(\varepsilon_h(\vec{v}^*))} \right) \snorm{\ojump{\vec{v}}} \dif s \nonumber \\
&\quad\quad\leq
2 \sum_{F \in \Gamma_I} \sum_{*\in\{+,-\}} \alpha^* w^* \norm{\sqrt{2 \mu} \varepsilon_h(\vec{v}^*)}_{L_2(\kappa^*)} \norm{\ojump{\vec{v}}}_{L_2(F)},
\end{align}
where
\begin{equation}
\alpha^* = C_\text{inv}(\kappa^*, F, p_{\kappa^*} - 2) \norm{2 \mu \vec{n}^*}_{L_\infty(F)} \norm{(2 \mu)^{-\frac{1}{2}}}_{L_\infty(\kappa^*)}.
\end{equation}
In order to proceed, following \eqref{eq_weights} and \eqref{eq_zeta}, we select
\begin{equation}
w^* = \frac{\zeta^*}{\zeta^+ + \zeta^-},
\quad
\zeta^* = \frac{1}{\delta \sqrt{m_{\kappa^*}} \alpha^*},
\end{equation}
respectively,
where $m_{\kappa^*}$ is the number of facets belonging to element $\kappa^*$, i.e., here $m_{\kappa^*}=3$, and
$\delta$ is a positive constant to be determined. Thereby, we deduce that
\begin{align}
&\left| 2 \int_{\Gamma_I} \avg{2 \mu \Pi_{L_2}(\varepsilon_h(\vec{v}))} \cdot \ojump{\vec{v}} \dif s \right| \nonumber \\
&\quad\quad\leq
2 \sum_{F \in \Gamma_I} \sum_{*\in\{+,-\}} \frac{1}{\delta \sqrt{m_{\kappa^*}}} \frac{1}{\zeta^+ + \zeta^-} \norm{\sqrt{2 \mu} \varepsilon_h(\vec{v}^*)}_{L_2(\kappa^*)} \norm{\ojump{\vec{v}}}_{L_2(F)} \nonumber \\
&\quad\quad\leq
2 \sum_{F\in\Gamma_I} \sum_{*\in\{+,-\}} \left( 
\frac{1}{2 \epsilon} \frac{1}{\delta^2 m_{\kappa^*}} \norm{\sqrt{2 \mu} \varepsilon_h(\vec{v}^*)}^2_{L_2(\kappa^*)}
+ \frac{\epsilon}{2} \beta \norm{\ojump{\vec{v}}}^2_{L_2(F)}
\right), \nonumber 
\end{align}
where $\epsilon > 0$ is a constant.
An analogous bound can also be established for faces on the boundary of $\Omega$. Hence, combining both bounds and
observing the inequality
\begin{equation}
\sum_{F \in \Gamma} \sum_{*\in\{+, -\}} \norm{z^*}^2_{L_2(\kappa^*)}
\leq \sum_{\kappa \in \mathcal{T}_h} m_\kappa \norm{z}^2_{L_2(\kappa)},
\end{equation}
gives
\begin{align}
\left| 2 \int_{\Gamma} \avg{2 \mu \Pi_{L_2}(\varepsilon_h(\vec{v}))} \cdot \ojump{\vec{v}} \dif s \right| 
\leq
\frac{1}{\delta^2 \epsilon} \sum_{\kappa \in \mathcal{T}^h} \norm{\sqrt{2 \mu} \varepsilon_h(\vec{v})}^2_{L_2(\kappa)}
+ 2 \epsilon \sum_{F \in \Gamma} \beta \norm{\ojump{\vec{v}}}^2_{L_2(F)},
\end{align}
from which we deduce that
\begin{equation}
\tilde{B}_\text{IP}(\vec{v}, \vec{v}) \geq 
\left(1 - \frac{1}{\delta^2\epsilon} \right) \norm{\sqrt{2 \mu} \varepsilon_h(\vec{v})}^2_{L_2(\Omega)}
+ \left(1 - 2 \epsilon \right) \norm{\sqrt{\beta} \ojump{\vec{v}}}^2_{L_2(\Gamma)}.
\end{equation}
To determine the unknown parameters $\delta$ and $\epsilon$ we set
\begin{equation}
\left(1 - \frac{1}{\delta^2\epsilon} \right) = \left(1 - 2 \epsilon \right)
\implies \epsilon = \frac{1}{\sqrt{2} \delta},
\end{equation}
which, provided $\delta > \sqrt{2}$, yields a coercive numerical scheme.

The proof of continuity follows based on employing a similar argument and hence the
details are omitted.
\end{proof}

\subsection{$C^0$-\gls{RIPG} convergence}
\label{sec:convergence}

In this section we now proceed to derive an $hp$-version a priori error bound for the
$C^0$-\gls{RIPG} formulation \eqref{eq:bilinear_fem}. For simplicity of
presentation, throughout this section we assume that $\mathcal{T}^h$ is a 
quasi-uniform triangular mesh with (global) mesh element size parameter $h$. For the
proceeding analysis we require an approximation result for the finite 
element space $V^{h,p}_0$. Given $v\in H^k(\Omega)$, $k\geq 3$, we assume 
there exists ${\mathcal I}: H^2(\Omega)\mapsto H^2(\Omega)\cap V^{h,p}_0$ such that
\begin{equation}
|v-{\mathcal I}v|_{H^s(\Omega)} \leq C \frac{h^{\mu-s}}{p^{k-s}} \|v\|_{H^k(\Omega)},
\label{eq:approx}
\end{equation}
where $\mu=\min(p+1,k)$ and $C$ is a positive constant independent of $h$, $p$,
$v$, and ${\mathcal I}v$. For triangular meshes in 2D this result was proved
in \citet{suri_1990} for $s=0,1,2$; see \citet{dong_mascotti_2023} for analogous approximation
results on tensor product meshes.

Equipped with \eqref{eq:approx} we now state the main result of this section.

\begin{theorem} \label{thm:apriori}
Given that ${\mathcal T}^h$ is a quasi-uniform triangular mesh, write
$\vec{u} = \nabla\times\phi$ and $\vec{u}_h = \nabla_h\times\phi_h$
to denote the solutions to \eqref{eq:stream_function_formulation} and 
\eqref{eq:bilinear_fem}, respectively. Then assuming
that $\phi \in H^k(\Omega)$, $k\geq 3$, we have that
$$
\dgnorm{\vec{u}-\vec{u}_h} \leq C \frac{h^{\mu-2}}{p^{k-2}} \|\phi\|_{H^k(\Omega)},
$$
where $\mu = \min(p+1,k)$ and $C$ is a positive constant, which is independent
of $h$ and $p$.
\end{theorem}

\begin{proof}
For $\vec{u} = \nabla\times\phi$ and $\vec{u}_h = \nabla_h\times\phi_h$ we recall Strang's lemma
\begin{align}
\dgnorm{\vec{u}-\vec{u}_h} 
&\leq \left(1+\frac{C_{\text{cont}}}{C_{\text{coer}}} \right) \inf_{\vec{v}_h = \nabla_h\times \psi_h, \psi_h \in V^{h,p}_0} \dgnorm{\vec{u}-\vec{v}_h}
\nonumber \\
&\hspace{1cm} + \frac{1}{C_{\text{coer}}}  \sup_{\vec{w}_h=\nabla_h\times \varphi_h, \varphi_h\in V^{h,p}_0\backslash\{0\}} 
\frac{| \tilde{B}_\text{IP}(\vec{u}, \vec{w}_h)-\ell (\vec{w}_h)|}{\dgnorm{\vec{w}_h}}.
\label{strang_inequality}
\end{align} 
Employing the approximation result \eqref{eq:approx}, together with the $H^2$-conformity of ${\mathcal I}$, the first term on the
right-hand side of \eqref{strang_inequality} can be bounded as follows
\begin{align*}
\inf_{\vec{v}_h = \nabla_h\times \psi_h, \psi_h \in V^{h,p}_0} \dgnorm{\vec{u}-\vec{v}_h}
& \leq \dgnorm{\nabla\times \phi-\nabla_h\times ({\mathcal I}\phi)} \\
& = \norm{\sqrt{2 \mu} \varepsilon_h(\nabla\times \phi - \nabla_h\times ({\mathcal I}\phi))}_{L_2(\Omega)} 
 \leq C \frac{h^{\mu-2}}{p^{k-2}} \|\phi\|_{H^k(\Omega)},
\end{align*}
as required.

For the second (consistency) term, noting that $\phi\in H^k(\Omega)$, $k\geq 3$, upon application of integration by parts
and the Cauchy-Schwarz inequality, we deduce that
\begin{align*}
\tilde{B}_\text{IP}(\vec{u}, \vec{w}_h)-\ell (\vec{w}_h)
& = \langle \avg{2 \mu (\varepsilon_h(\vec{u}) - \Pi_{L_2} \varepsilon_h(\vec{u}))}, \ojump{\vec{w}_h} \rangle_{\Gamma} \\
& \leq \| \beta^{-1/2} \avg{2 \mu (\varepsilon_h(\vec{u}) - \Pi_{L_2} \varepsilon_h(\vec{u}))} \|_{L_2(\Gamma)}
       \| \sqrt{\beta} \ojump{\vec{w}_h} \|_{L_2(\Gamma)}.
\end{align*}
We now recall the following approximation result from \citet{chernov_2012} for the $L_2$-projector. With a slight abuse of notation, 
we also write $\Pi_{L_2}:L_2(\kappa) \mapsto \mathcal{P}_p(\kappa)$  to denote the elementwise (scalar-valued) $L_2$-projector. 
With this notation, given a triangular element
$\kappa \in {\mathcal T}^h$, let $F\subset \partial \kappa$ denote one of its faces. Then, for $v\in H^k(\kappa)$, $k\geq 1$, the
following bound holds
\begin{equation}
\| v - \Pi_{L_2} v \|_{L_2(F)} \leq C \frac{h^{\mu-1/2}}{p^{k-1/2}} \|v\|_{H^k(\kappa)},
\label{eq:approx2}
\end{equation}
where $\mu=\min(p+1,k)$ and $C$ is a positive constant independent of $h$, $p$,
$v$, and $\Pi_{L_2} v$. 

Equipped with \eqref{eq:approx2}, the definition of $\beta$, cf.~\eqref{eq_weights} and the DG norm \eqref{eq:ipnorm}, we get
\begin{align*}
\tilde{B}_\text{IP}(\vec{u}, \vec{w}_h)-\ell (\vec{w}_h)
\leq C \frac{h^{\mu-2}}{p^{k-3/2}} \|\phi\|_{H^k(\Omega)} \dgnorm{w_h}.
\end{align*}
Collecting the above bounds gives the desired result.
\end{proof}

\begin{remark}
We remark that the bound derived in Theorem~\ref{thm:apriori} is optimal in both the mesh element size $h$ 
and the polynomial degree $p$; this is in agreement to the analogous bound derived for the
standard $C^0$-\gls{IPG} scheme in \citet{dong_mascotti_2023} for the Dirichlet problem. In the
case when inhomogeneous Dirichlet boundary conditions are employed, then as in the case of
second-order linear elliptic partial differential equations $p$-optimality is no longer possible,
cf.~\cite{Georgoulis_et_al_2009}.
\end{remark}

\section{Numerical experiments}
\label{sec:experiments}

In this section we present a series of numerical experiments to investigate the
practical performance of the proposed $C^0$-\gls{RIPG} scheme. We note that
all of the computational examples have been implemented using
the components of the FEniCS project~\citep{fenics:book}. We highlight the
\gls{UFL}~\citep{ufl} in particular as it facilitates the straightforward
specification of the verbose facet terms arising in the $C^0$-\gls{RIPG}
formulation \eqref{eq:bilinear_fem}. Initial prototypes of the
$C^0$-\gls{IPG} and \gls{SIPG} formulations were developed with the principles
outlined in \citet{nate2018}. The code used to generate the results presented
here is available in \citet{nate2023C0RIPGcode}.

Our experiments are constructed in two settings to test the numerical scheme.
Firstly, we present a numerical example with a known analytical solution in
order to validate the optimality of the a priori error bound derived in
\Cref{thm:apriori}. Secondly, we study the performance of the proposed
$C^0$-\gls{RIPG} scheme in the practical setting of reproducing benchmarks in
mantle convection cell models. We examine the benefits of the $C^0$-\gls{RIPG}
scheme's pointwise divergence free velocity approximation when coupled with
advection of a scalar (temperature) field, in addition to its performance with
viscosity models which are composed of variations over many orders of
magnitude.

\subsection{Manufactured solution}
\label{sec:mms_solution}

\begin{figure}[t!]
\centering
\includegraphics[width=0.8\linewidth]{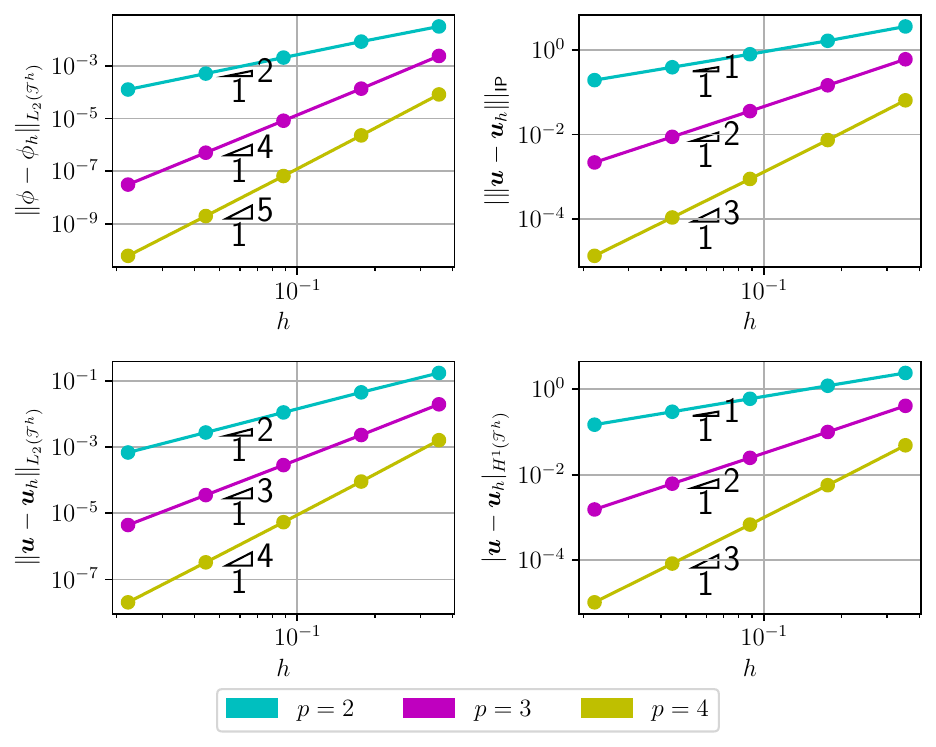}
\caption{Measured error norms of the manufactured solution problem.}
\label{fig:mms_convergence}
\end{figure}

In this section, we let $\Omega = (-1, 1)^2$ be a square. We subdivide $\Omega$ into a hierarchy
of meshes composed of $N \times N$ shape regular quadrilaterals each
bisected into triangle elements, where $N \in \{8, 16, 32, 64, 128\}$ is the number of
quadrilaterals dividing each orthogonal direction of $\Omega$. We select the
analytical solution and viscosity to be
\begin{equation}
\phi = \pi^{-1} \sin(\pi x) \sin(\pi y) \text{ and } \mu = 1 + \sin^2(\pi x) \sin^2(\pi y),
\end{equation}
respectively, such that
\begin{equation}
\vec{u} = \begin{pmatrix}
\sin(\pi x) \cos(\pi y) \\
- \cos(\pi x) \sin(\pi y) \\
\end{pmatrix},
\end{equation}
which then determines $\vec{f}$ according to \cref{eq:momentum}. We compute
approximations of $\phi$ using the $C^0$-\gls{RIPG} formulation and examine
the rate at which the numerical approximation converges to the analytical solution.
Furthermore, we examine the influence of the parameter $\delta$ on the stability
of the numerical scheme. Here, the approximation error is measured in the DG-norm
\eqref{eq:ipnorm}, as well as the following norms of interest:
\begin{align}
\norm{\phi - \phi_h}^2_{L_2(\mathcal{T}^h)} 
&= \sum_{\kappa \in \mathcal{T}^h} \int_\kappa (\phi - \phi_h)^2 \mathrm{d} \vec{x}, \\
\norm{\vec{u} - \vec{u}_h}^2_{L_2(\mathcal{T}^h)} 
&= \sum_{\kappa \in \mathcal{T}^h} \int_\kappa (\vec{u} - \nabla \times \phi_h)^2 \mathrm{d} \vec{x}, \\
\snorm{\vec{u} - \vec{u}_h}^2_{H^1(\mathcal{T}^h)}
&= \sum_{\kappa \in \mathcal{T}^h} \int_\kappa \left(\nabla \vec{u} - \nabla(\nabla \times \phi_h)\right)^2 \mathrm{d} \vec{x}.
\end{align}
In \Cref{fig:mms_convergence} we plot the error measured in the above norms against the mesh element size $h$ on the aforementioned sequence of
uniform (structured) triangular meshes for $p=2,3,4$. Here, we observe that for each fixed $p$, the DG norm of the error tends to zero at the optimal
rate of ${\mathcal O}(h^{p-1})$ as $h$ tends to zero; this is in full agreement with the predicted rate given in \Cref{thm:apriori}. 
Analogous rates are observed for $\snorm{\vec{u} - \vec{u}_h}^2_{H^1(\mathcal{T}^h)}$ as expected. The $L_2(\Omega)$ norm of the error
in the approximation to the velocity is observed to behave like ${\mathcal O}(h^{p})$ as $h$ tends to zero, which are again expected. However,
we observe that $\norm{\phi - \phi_h}^2_{L_2(\mathcal{T}^h)} = {\mathcal O}(h^{2})$ for $p=2$, while 
$\norm{\phi - \phi_h}^2_{L_2(\mathcal{T}^h)} = {\mathcal O}(h^{p+1})$ for $p=3,4$, as $h$ tends to zero, which indicates suboptimality in the
approximation of the stream function when the error is measured in $L_2(\Omega)$ norm and the lowest order approximation is employed.

\begin{figure}[t!]
\centering
\includegraphics[width=0.8\linewidth]{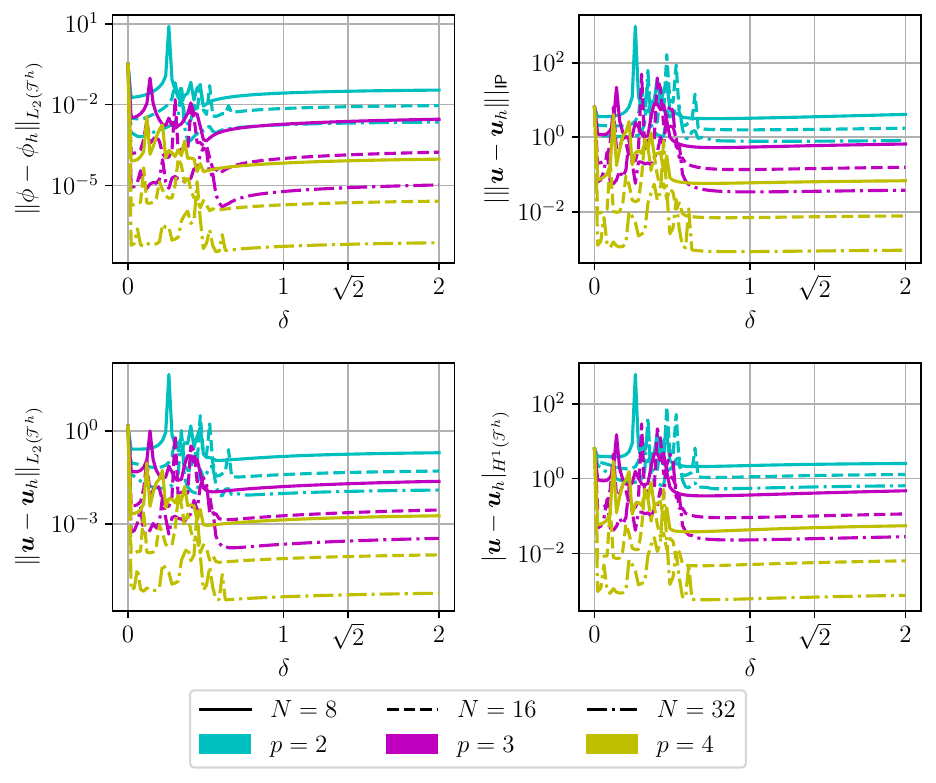}
\caption{Measured error norms of the manufactured solution problem with varying
parameter $\delta$. Here we ensure the $C^0$-\gls{RIPG} scheme is stable with
$\delta > \sqrt{2}$ as predicted in \Cref{sec:coercivity}.}
\label{fig:mms_sigmatest}
\end{figure}

In \Cref{fig:mms_sigmatest} we study the influence of the parameter $\delta$ on the stability of the $C^0$-\gls{RIPG} scheme and the approximation error.
As indicated in \Cref{sec:coercivity} we observe that setting the penalty parameter $\delta > \sqrt{2}$ guarantees stability of the
underlying method. We also see that smaller values of $\delta$ may be employed in practice which lead to a slight improvement
of the error measured in the above norms. However, as expected if $\delta$ is reduced too far, then the stability of the $C^0$-\gls{RIPG} scheme is
no longer guaranteed.


{}

\subsection{Buoyancy-driven flow}
\label{sec:tosi_problem}

To examine the practicality of the $C^0$-\gls{RIPG} method we reproduce
numerical benchmarks in the context of geophysical flow driven by thermal
buoyancy. These benchmarks further employ temperature and strain-rate
dependent viscosity models as exhibited in the works \citet{blankenbach1989}
and \citet{tosi2015}. These benchmarks require that $\Omega = (0, 1)^2$ is the
unit square, in which we seek the velocity and pressure which satisfy
\cref{eq:momentum,eq:mass} subject to the free slip boundary conditions
\begin{alignat}{2}
\vec{u} \cdot \vec{n} &= 0 &&\text{ on } \partial\Omega, \\
\left(2 \mu \varepsilon(\vec{u}) \cdot \vec{n} - P \vec{n} \right) \cdot \vec{\tau} &= 0 &&\text{ on } \partial\Omega.
\end{alignat}
For thermally driven buoyancy, we set
\begin{equation}
\vec{f} = \mathrm{Ra} \, T  \,\hat{\vec{y}},
\end{equation}
with viscosity model
\begin{align}
\mu &= 
\begin{cases}
\mu_\text{lin} & \sigma_Y = 0, \\
2 (\mu_\text{lin}^{-1} + \mu_\text{plast}^{-1})^{-1} & \sigma_Y > 0,
\end{cases} \\
\mu_\text{lin} &= \exp\left(-\log (\Delta \mu_T) T + \log(\Delta \mu_z) z\right), \\
\mu_\text{plast} &= \num{e-3} + \frac{\sigma_Y}{\sqrt{\varepsilon(\vec{u}) : \varepsilon(\vec{u})}}.
\end{align}
Here, $z = 1 - y$ is a measure of depth, $\mathrm{Ra}$ is the
constant Rayleigh number and $\Delta \mu_T$, $\Delta \mu_z$ and $\sigma_Y$ are
constant viscosity model parameters to be defined in each benchmark case. The temperature field 
$T:\Omega \mapsto \mathbb{R}^+ \cup \{ 0 \}$
satisfies the following advection--diffusion problem
\begin{alignat}{2}
\vec{u} \cdot \nabla T - \nabla^2 T &= 0 \quad &&\text{ in } \Omega, 
\label{eq:heat} \\
T &= 0 \quad &&\text{ on } [0, 1] \times \{ 0 \}, 
\label{eq:heat_bc0} \\
T &= 1 \quad &&\text{ on } [0, 1] \times \{ 1 \}, 
\label{eq:heat_bc1} \\
\nabla T \cdot \vec{n} &= 0 \quad &&\text{ on } \{0, 1\} \times [0, 1].
\label{eq:heat_bc2}
\end{alignat}
We discretize \cref{eq:heat,eq:heat_bc0,eq:heat_bc1,eq:heat_bc2} using standard
$C^0$~finite elements such that we seek $(\phi_h, T_h) \in
V^{h,p}_0 \times S^{h,p}$ such that \cref{eq:bilinear_fem} holds simultaneously with
\begin{equation}
(\vec{u}_h \cdot \nabla T_h, s_h) + (\nabla T_h, \nabla s_h) = 0
\end{equation}
for all $(\psi_h, s_h) \in V^{h,p}_0 \times  S^{h,p}_0$, where
$S^{h,p} = \{ v \in H^1(\Omega) : \left. v \right|_{\kappa} \in \mathcal{P}_p(\kappa) \:\: \forall \kappa \in \mathcal{T}^h, \:\: \left. v \right|_{[0, 1] \times \{0\}} = 0, \:\: \left. v \right|_{[0, 1] \times \{1\}} = 1\}$ and $S^{h,p}_0 = \{ v \in H^1(\Omega) : \left. v \right|_{\kappa} \in \mathcal{P}_p(\kappa) \:\: \forall \kappa \in \mathcal{T}^h, \:\: \left. v \right|_{[0, 1] \times \{0, 1\}} = 0 \}$.

\begin{table}[t!]
\centering
\begin{tabular}{lrrrrrr}
\toprule
Case & $\mathrm{Ra}$ & $\Delta \mu_T$ & $\Delta \mu_z$ & $\sigma_Y$ & $\mathrm{Nu}_\mathrm{ref}$ & $u_\mathrm{rms,ref}$ \\ \midrule
BB1a & \num{e4} & \num{1} & \num{1} & \num{0} & \num{4.884409} & \num{42.864947} \\
BB2a & \num{e4} & \num{e4} & \num{1} & \num{0} & \num{10.065899} & \num{480.433425}\\
T2   & \num{e2} & \num{e5} & \num{1} & \num{1} & \num{8.559459} & \num{140.775535} \\ 
T4   & \num{e2} & \num{e5} & \num{10} & \num{1} & \num{6.615419} & \num{79.088809} \\
\bottomrule                 
\end{tabular}
\caption{Benchmark cases exhibited in \protect\citet{blankenbach1989} and
\protect\citet{tosi2015} (BB and T prefixes, respectively) and corresponding reference values selected from
those works.}
\label{tab:benchmark_cases}
\end{table}

\begin{table}[t!]
\centering
\begin{tabular}{llrrrrrr}
\toprule
Case & $p$ & $N$ & $\mathrm{Nu}$ & $u_\mathrm{rms}$ & $\langle W \rangle$ & $\langle \Phi \rangle / \mathrm{Ra}$ & $\Delta$ \\
\midrule
\multirow[t]{10}{*}{BB1a} & \multirow[t]{5}{*}{2} & \num{8} & \num{5.159790} & \num{41.877692} & \num{3.812435} & \num{3.721304} & \num{2.390357e-02} \\
 &  & \num{16} & \num{5.015557} & \num{42.613520} & \num{3.866134} & \num{3.841287} & \num{6.426728e-03} \\
 &  & \num{32} & \num{4.923947} & \num{42.801576} & \num{3.879774} & \num{3.873409} & \num{1.640653e-03} \\
 &  & \num{64} & \num{4.894819} & \num{42.849068} & \num{3.883246} & \num{3.881644} & \num{4.124164e-04} \\
 &  & \num{128} & \num{4.887047} & \num{42.860973} & \num{3.884118} & \num{3.883717} & \num{1.032472e-04} \\
 \cmidrule{2-8} & \multirow[t]{5}{*}{3} & \num{8} & \num{4.997386} & \num{42.859417} & \num{3.883614} & \num{3.894593} & \num{2.819064e-03} \\
 &  & \num{16} & \num{4.894932} & \num{42.864517} & \num{3.884342} & \num{3.887360} & \num{7.764650e-04} \\
 &  & \num{32} & \num{4.885120} & \num{42.864918} & \num{3.884405} & \num{3.885177} & \num{1.986228e-04} \\
 &  & \num{64} & \num{4.884454} & \num{42.864943} & \num{3.884409} & \num{3.884603} & \num{4.993627e-05} \\
 &  & \num{128} & \num{4.884412} & \num{42.864945} & \num{3.884409} & \num{3.884458} & \num{1.250158e-05} \\
 \midrule \multirow[t]{10}{*}{BB2a} & \multirow[t]{5}{*}{2} & \num{8} & \num{11.569494} & \num{494.791616} & \num{10.639110} & \num{9.866476} & \num{7.262202e-02} \\
 &  & \num{16} & \num{10.731001} & \num{487.372110} & \num{9.151660} & \num{8.917624} & \num{2.557305e-02} \\
 &  & \num{32} & \num{10.306838} & \num{483.106619} & \num{9.101327} & \num{9.018523} & \num{9.098057e-03} \\
 &  & \num{64} & \num{10.129296} & \num{481.246438} & \num{9.072895} & \num{9.046327} & \num{2.928342e-03} \\
 &  & \num{128} & \num{10.081744} & \num{480.619011} & \num{9.067279} & \num{9.059942} & \num{8.091935e-04} \\
 \cmidrule{2-8} & \multirow[t]{5}{*}{3} & \num{8} & \num{11.039987} & \num{481.181654} & \num{9.217020} & \num{9.258907} & \num{4.524017e-03} \\
 &  & \num{16} & \num{10.225818} & \num{481.626650} & \num{9.107274} & \num{9.123698} & \num{1.800123e-03} \\
 &  & \num{32} & \num{10.073773} & \num{480.583438} & \num{9.064849} & \num{9.072503} & \num{8.436236e-04} \\
 &  & \num{64} & \num{10.066181} & \num{480.403057} & \num{9.065588} & \num{9.067885} & \num{2.533984e-04} \\
 &  & \num{128} & \num{10.065910} & \num{480.427513} & \num{9.065872} & \num{9.066485} & \num{6.764788e-05} \\
 \midrule \multirow[t]{10}{*}{T2} & \multirow[t]{5}{*}{2} & \num{8} & \num{8.670484} & \num{130.718863} & \num{7.175272} & \num{7.024771} & \num{2.097492e-02} \\
 &  & \num{16} & \num{8.852719} & \num{135.810284} & \num{7.405343} & \num{7.330546} & \num{1.010052e-02} \\
 &  & \num{32} & \num{8.677987} & \num{139.186405} & \num{7.514772} & \num{7.477494} & \num{4.960578e-03} \\
 &  & \num{64} & \num{8.587873} & \num{140.198404} & \num{7.542252} & \num{7.528339} & \num{1.844703e-03} \\
 &  & \num{128} & \num{8.565231} & \num{140.580896} & \num{7.553466} & \num{7.548985} & \num{5.933286e-04} \\
 \cmidrule{2-8} & \multirow[t]{5}{*}{3} & \num{8} & \num{9.044719} & \num{137.032238} & \num{7.479251} & \num{7.480751} & \num{2.004897e-04} \\
 &  & \num{16} & \num{8.622321} & \num{140.017268} & \num{7.533352} & \num{7.536932} & \num{4.750446e-04} \\
 &  & \num{32} & \num{8.557657} & \num{140.503070} & \num{7.550360} & \num{7.552252} & \num{2.506208e-04} \\
 &  & \num{64} & \num{8.557513} & \num{140.702836} & \num{7.557034} & \num{7.557831} & \num{1.054292e-04} \\
 &  & \num{128} & \num{8.559291} & \num{140.769772} & \num{7.559261} & \num{7.559543} & \num{3.730833e-05} \\
 \midrule \multirow[t]{10}{*}{T4} & \multirow[t]{5}{*}{2} & \num{8} & \num{6.964560} & \num{76.225863} & \num{5.486364} & \num{5.347143} & \num{2.537569e-02} \\
 &  & \num{16} & \num{6.879954} & \num{78.177425} & \num{5.583371} & \num{5.521602} & \num{1.106295e-02} \\
 &  & \num{32} & \num{6.702668} & \num{78.807864} & \num{5.604996} & \num{5.583245} & \num{3.880586e-03} \\
 &  & \num{64} & \num{6.637644} & \num{78.990505} & \num{5.611148} & \num{5.604410} & \num{1.200818e-03} \\
 &  & \num{128} & \num{6.620866} & \num{79.060030} & \num{5.614094} & \num{5.612224} & \num{3.329760e-04} \\
 \cmidrule{2-8} & \multirow[t]{5}{*}{3} & \num{8} & \num{6.931491} & \num{78.787124} & \num{5.613935} & \num{5.617167} & \num{5.754950e-04} \\
 &  & \num{16} & \num{6.651302} & \num{79.002440} & \num{5.610692} & \num{5.614473} & \num{6.734658e-04} \\
 &  & \num{32} & \num{6.615929} & \num{79.045743} & \num{5.612879} & \num{5.614332} & \num{2.588058e-04} \\
 &  & \num{64} & \num{6.615222} & \num{79.082175} & \num{5.615025} & \num{5.615510} & \num{8.630228e-05} \\
 &  & \num{128} & \num{6.615397} & \num{79.088264} & \num{5.615384} & \num{5.615517} & \num{2.360951e-05} \\
\bottomrule
\end{tabular}
\caption{Computed functional measurements from the numerical benchmark cases as computed with the $C^0$-\gls{RIPG} scheme.}
\label{tab:benchmark_functionals}
\end{table}

\begin{figure}[t!]
\centering
\includegraphics[width=0.8\linewidth]{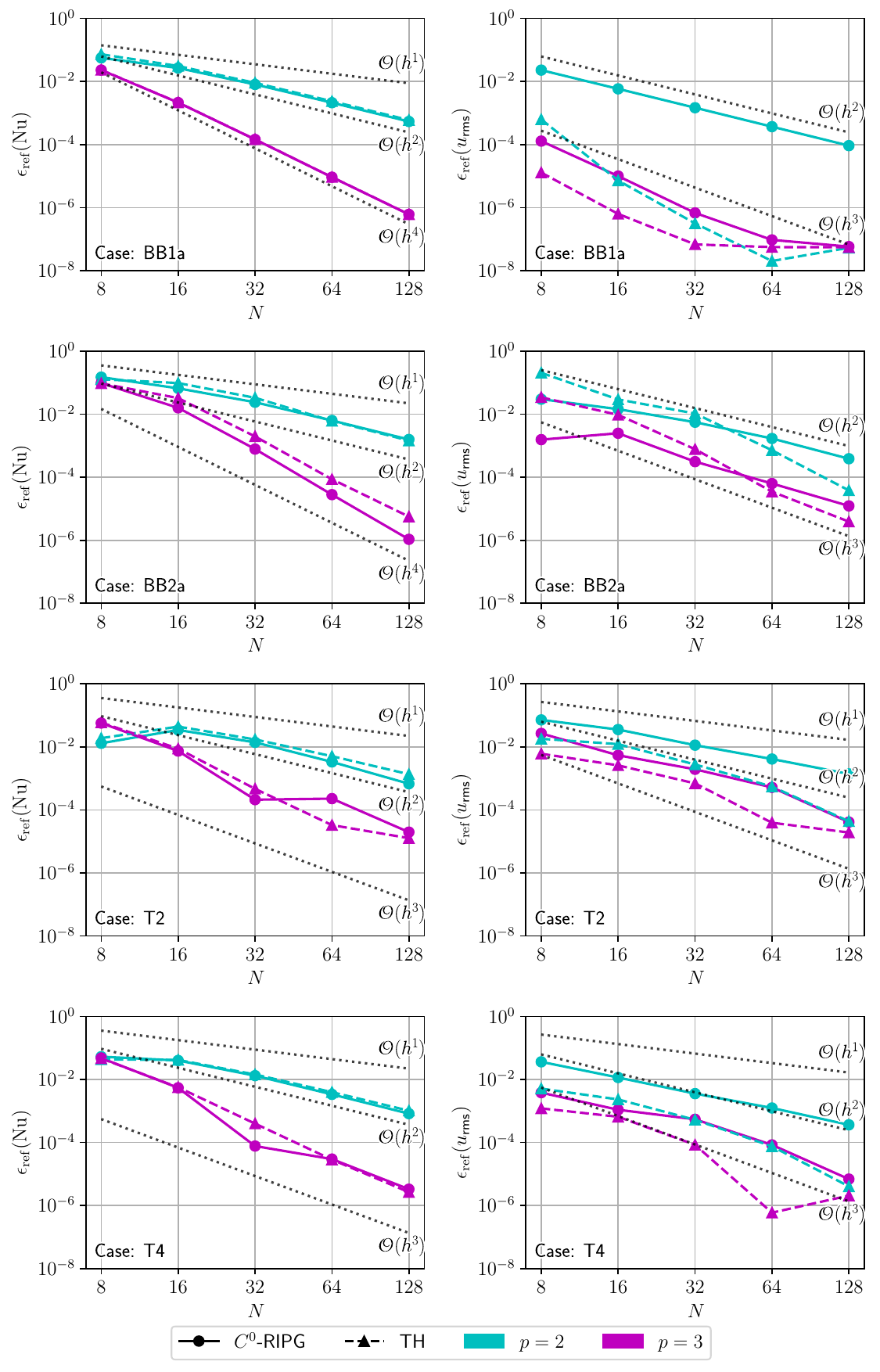}
\caption{
    Measured functionals' errors relative to reference values selected from
    benchmark problems exhibited in \protect\citet{blankenbach1989,tosi2015}.
    We further compare the $C^0$-\gls{RIPG} formulation with a standard
    \gls{TH} discretization.
    }
\label{fig:blankenbach_convergence}
\end{figure}

\begin{figure}[t!]
\centering
\includegraphics[width=0.7\linewidth]{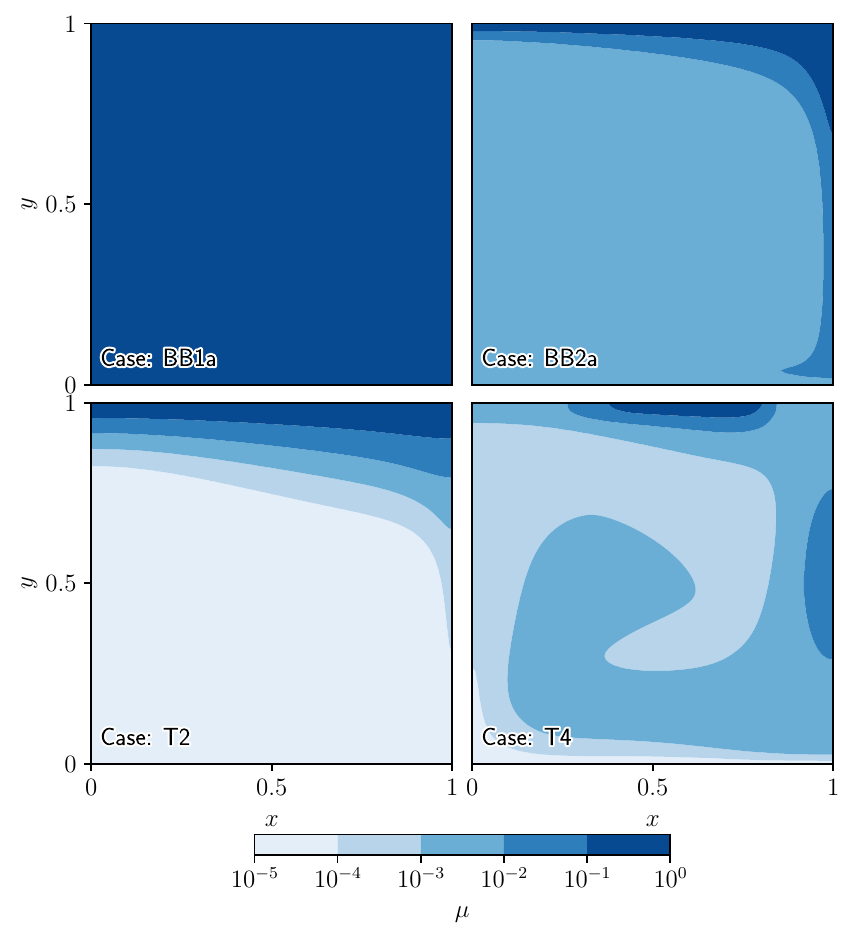}
\caption{
    Viscosity fields as computed from the $C^0$-\gls{RIPG} discretization with
    $p=3$ of the benchmark problems tabulated in \Cref{tab:benchmark_cases}.
    These fields are interpolated in the space $V^{h,p=1}_\text{DG}$ for
    visualization. }
\label{fig:blankenbach_viscosity}
\end{figure}

In \Cref{tab:benchmark_cases} we define a number of benchmark cases with
corresponding Rayleigh numbers and viscosity models. The functionals of
interest measured for comparison with the corresponding benchmark reports are
as follows
%
%
\begin{alignat}{6}
\phantom{\left\langle T \right\rangle} &\phantom{= \sum_{\kappa \in {\mathcal{T}^h}} \int_\kappa T_h \, \mathrm{d}\vec{x},} &\quad
\mathrm{Nu} &= \int_0^1 \left.(\nabla T_h \cdot \vec{n})\right|_{y=1} \, \mathrm{d} x, &\quad
u_\mathrm{rms} &= \sqrt{\int_\Omega \vec{u}^2_h \, \mathrm{d} \vec{x}}, \nonumber \\
\left\langle W \right\rangle &= \sqrt{\int_\Omega T_h \vec{u}_h \cdot \hat{\vec{y}} \, \mathrm{d}\vec{x}}, &\quad
\left\langle \Phi\right\rangle &= \sqrt{\int_\Omega 2 \mu \varepsilon_h(\vec{u}_h) : \varepsilon_h(\vec{u}_h) \, \mathrm{d} \vec{x}}, &\quad
\Delta &= 
\frac{\left\vert \left\langle W \right\rangle - \frac{\left\langle \Phi\right\rangle}{\mathrm{Ra}} \right\vert}
{\max \left( \left\langle W \right\rangle, \frac{\left\langle \Phi\right\rangle}{\mathrm{Ra}} \right)}, \nonumber
\end{alignat}
where 
$\mathrm{Nu}$ is the Nusselt number at the top boundary,
$u_\text{rms}$ is the root mean square speed, $\left\langle W \right\rangle$ and $\left\langle \Phi \right\rangle$ are
the average rates of work done against gravity and viscous dissipation,
respectively, and $\Delta$ is a measure of thermal energy conservation.
In \Cref{tab:benchmark_functionals} we tabulate the computed functional
values from our implementation on a sequence of uniform triangular meshes for $p=2,3$.

In \Cref{fig:blankenbach_convergence} we compute errors in $\mathrm{Nu}$ and $u_\mathrm{rms}$
relative to the given reference values for
the benchmark cases shown in \Cref{tab:benchmark_cases}. Here, the relative error
measurement of computed quantity $\chi$ compared with reference value
$\chi_\mathrm{ref}$ is given by
\begin{equation}
\epsilon_\mathrm{ref}(\chi) = \frac{\vert \chi - \chi_\mathrm{ref} \vert}{\chi_\mathrm{ref}}.
\end{equation}

For the case BB1a (isoviscous model) we see that as the mesh is uniformly refined, 
the computed Nusselt number and root-mean-square velocity approximations, as compared with reference values,
converge to zero at the rates ${\mathcal O}(h^{2(p-1)})$ and ${\mathcal O}(h^{p})$ as $h$ tends to zero, respectively, for $p=2,3$. 
However, we observe that when nonlinear viscosity models are employed, as in cases BB2a, T2 and T4, the order
of convergence of these quantities of interest may be reduced; this is particularly evident for
the `highly nonlinear' cases T2 and T4 (cf. \Cref{fig:blankenbach_viscosity}).

Additionally, we compare the performance of the proposed $C^0$-\gls{RIPG} method with the
\gls{TH} discretization scheme which highlights the impact of employing an exactly
divergence free velocity field approximation. When examining equivalent
polynomial degrees in the underlying finite element spaces, we see that the
$C^0$-\gls{RIPG} formulation yields a more precise Nusselt number
approximation for fewer than half the number of degrees of freedom due to its
solenoidal velocity field. This result is attained despite the fact that the
root-mean-square velocity approximation of the $C^0$-\gls{RIPG} scheme is less
accurate than that the corresponding quantity computed for the \gls{TH} method
for an equivalent polynomial degree. For the $C^0$-\gls{RIPG} method to
yield a root-mean-square velocity approximation roughly equivalent to
the \gls{TH} method a polynomial degree one order higher should be employed.

Finally, we examine the balance of energy encapsulated by the functional
$\Delta$ as tabulated in \Cref{tab:benchmark_functionals}. Here, we see that energy
is not exactly conserved. However, the values of $\Delta$ computed using the
$C^0$-\gls{RIPG} scheme are comparable with the \gls{TH} scheme results,
$\Delta<0.01\%$, reported in \citet{tosi2015}.

\section{Concluding remarks}
\label{sec:conclusion}

In this article we have presented the $C^0$-\gls{RIPG} formulation of the stream function
formulation of the Stokes system with varying viscosity. We have shown that this
scheme is stable provided the interior penalty parameter is selected so that 
$\delta > \sqrt{2}$. Moreover, our analysis and numerical experiments demonstrate that
the scheme converges optimally with respect to uniform mesh refinement when the error 
is measured in an appropriate DG norm. Furthermore, the discretization provides an
exactly divergence free approximation of the velocity which, in the context of
numerical benchmarks of mantle convection simulations, yields a more precise
approximation of the advected temperature field. Our implementation is
available in \citet{nate2023C0RIPGcode} for use with the FEniCSx library.

\section{Acknowledgments}
\label{sec:acknowledgments}

This research was partly supported by NSF-EAR grant 2021027. N. Sime
gratefully acknowledges the support of Carnegie Institution for Science
President's Fellowship.

\appendix
\section{$C^0$-\gls{RIPG} formulation derivation}
\label{sec:appdx_derivation}

In this appendix we outline the derivation of the proposed $C^0$-\gls{RIPG}
scheme \eqref{eq:bilinear_fem}; to this end we first develop the \gls{SIPG}
\gls{DG} discretization of the problem in
\cref{eq:momentum,eq:mass,eq:bcs-1,eq:bcs-2}, then consider the restriction to
the $H^1$-conforming FE space $V^{h,p}$.

\subsection{Flux formulation}
%
For completeness, we consider the following generalized stream function formulation of \cref{eq:momentum} with corresponding boundary
conditions:
\begin{align}
\nabla \times \left( -\nabla \cdot \left(2 \mu \varepsilon(\nabla \times \phi) - P \mathbb{I} \right) \right)
&=
\nabla \times \vec{f} && \text{in } \Omega, \label{eq:sf_eqn}\\
\phi &= \phi_D && \text{on } \partial\Omega_{\phi,D} \\
\nabla \times \phi &= \vec{u}_D && \text{on } \partial\Omega_{D}, \\
\left(2 \mu \varepsilon(\nabla \times \phi) - P \mathbb{I} \right) \cdot \vec{n} &= \vec{g}_{N} && \text{on } \partial\Omega_{N}, \\
\left(-\nabla \cdot \left(2 \mu \varepsilon(\nabla \times \phi) - P \mathbb{I} \right)\right) \times \vec{n} &= \vec{f} \times \vec{n} && \text{on } \partial\Omega_{\phi,N}.
\end{align}
Here, $\mathbb{I}$ is the identity tensor and the exterior boundary is split into
components for the stream function and velocity boundary conditions such that
$\partial\Omega = \partial\Omega_{\phi,D} \cup \partial\Omega_{\phi,N}$ and
$\partial \Omega = \partial\Omega_{D} \cup \partial\Omega_{N}$
where the boundary components do not overlap, $\partial\Omega_{\phi,D} \cap
\partial\Omega_{\phi,N} = \emptyset$ and $\partial\Omega_{D} \cap
\partial\Omega_{N} = \emptyset$. We define the rank 4 tensor $G =
\partial (2 \mu \varepsilon(\vec{u})) / \partial (\nabla
\vec{u})$; given $\sigma \in \mathbb{R}^{2\times
2}$, its product and transpose product are defined, respectively, by
\begin{equation}
G \sigma = G_{ijkl} \sigma_{ij} \text{ and } G^\top \sigma = G_{ijkl} \sigma_{kl},
\end{equation}
such that $2 \mu \varepsilon(\vec{u}) = G \nabla \vec{u}$. 
Let $F_1 : \Omega \mapsto \mathbb{R}^2$, $F_2 : \Omega \mapsto
\mathbb{R}^{2 \times 2}$ and $F_3 : \Omega \mapsto \mathbb{R}^2$ such that
we recast \cref{eq:sf_eqn} in terms of four first order equations
\begin{equation}
\nabla \times \vec{f} = \nabla \times F_1, \quad
F_1 = -\nabla \cdot F_2, \quad
F_2 = G \nabla F_3 - P \mathbb{I}, \quad
F_3 = \nabla \times \phi,
\label{eq:sf_eqn_firstorder}
\end{equation}
where it is evident that
\begin{equation}
F_1(\phi) = -\nabla \cdot \left(2 \mu \varepsilon(\nabla \times \phi) - P \mathbb{I} \right) ~~ \mbox{ and } ~~
F_2(\phi) =  2 \mu \varepsilon (\nabla \times \phi) - P \mathbb{I}.
\end{equation}
We multiply each equation in \eqref{eq:sf_eqn_firstorder} by $v_1 : \Omega \mapsto
\mathbb{R}$, $\vec{v}_2 : \Omega \mapsto
\mathbb{R}^{2}$, $v_3 : \Omega \mapsto \mathbb{R}^{2\times 2}$ and
$\vec{v}_4 : \Omega \mapsto \mathbb{R}^2$, respectively, and integrate over an element
$\kappa \in \mathcal{T}^h$ to give
\begin{align}
(\nabla \times \vec{f}, v_1)_\kappa &= (\nabla \times F_1, v_1)_\kappa, \label{eq:tfmult_1} \\
(F_1, \vec{v}_2)_\kappa &= (-\nabla \cdot F_2, \vec{v}_2)_\kappa, \label{eq:tfmult_2} \\
(F_2, v_3)_\kappa &= (G \nabla F_3, v_3)_\kappa - (P \mathbb{I}, v_3)_\kappa, \label{eq:tfmult_3} \\
(F_3, \vec{v}_4)_\kappa &= (\nabla \times \phi, \vec{v}_4)_\kappa. \label{eq:tfmult_4}
\end{align}
Integrating by parts \cref{eq:tfmult_1,eq:tfmult_2,eq:tfmult_3} once, and
given the lack of prescribed boundary data, \cref{eq:tfmult_4}
twice, we deduce that
\begin{align}
(\nabla \times \vec{f}, v_1)_\kappa &= (F_1, \nabla \times v_1)_\kappa
- \langle \widehat{F_1}, \vec{n} \times v_1 \rangle_{\partial\kappa}, \label{eq:ibp_1} \\
(F_1, \vec{v}_2)_\kappa &= (F_2, \nabla \vec{v}_2)_\kappa 
- \langle \widehat{F_2}, \vec{v}_2 \otimes \vec{n} \rangle_{\partial\kappa}, \label{eq:ibp_2}\\
(F_2, v_3)_\kappa &= -(F_3, \nabla \cdot (G^\top v_3))_\kappa 
+ \langle \widehat{F_3} \otimes \vec{n}, G^\top v_3 \rangle_{\partial\kappa} - (P \mathbb{I}, v_3)_\kappa, \label{eq:ibp_3}\\
(F_3, \vec{v}_4)_\kappa &= (\nabla \times \phi, \vec{v}_4)_\kappa 
- \langle \widehat{\phi} - \phi, \vec{n} \times \vec{v}_4 \rangle_{\partial\kappa}, \label{eq:ibp_4}
\end{align}
where $\widehat{(\cdot)}|_{\partial\kappa}$ indicates a consistent and
conservative flux approximation. We now proceed to eliminate the additional auxiliary variables 
introduced in order to derive the so-called flux formulation. To this end, we first select
$\vec{v}_4 = \nabla \cdot (G^\top v_3)$; inserting \cref{eq:ibp_4} into
\cref{eq:ibp_3} yields
\begin{equation}
(F_2, v_3)_\kappa = -(\nabla \times \phi, \nabla \cdot (G^\top v_3))_\kappa 
+ \langle \widehat{\phi} - \phi, \vec{n} \times (\nabla \cdot (G^\top v_3)) \rangle_{\partial\kappa}
+ \langle \widehat{F}_3 \otimes \vec{n}, G^\top v_3 \rangle_{\partial\kappa} - (P \mathbb{I}, v_3)_\kappa. \label{eq:midway_1}
\end{equation}
Given a lack of prescribed boundary information regarding $F_3$ we integrate the first term in 
\cref{eq:midway_1} by parts again giving
\begin{equation}
(F_2, v_3)_\kappa = (\nabla(\nabla \times \phi_h), G^\top v_3)_\kappa
+ \langle \widehat{\phi} - \phi, \vec{n} \times (\nabla \cdot (G^\top v_3)) \rangle_{\partial\kappa}
+ \langle (\widehat{F}_3 - \nabla \times \phi) \otimes \vec{n}, G^\top v_3 \rangle_{\partial\kappa}
- (P \mathbb{I}, v_3)_\kappa. \label{eq:midway_2}
\end{equation}
Let $\vec{v}_2 = \nabla \times v_1$, then inserting \cref{eq:ibp_2} into \cref{eq:ibp_1} yields
\begin{equation}
(\nabla \times \vec{f}, v_1)_\kappa = (F_2, \nabla (\nabla \times v_1))_\kappa 
- \langle \widehat{F}_2, (\nabla \times v_1) \otimes \vec{n} \rangle_{\partial\kappa}
- \langle \widehat{F}_1, \vec{n} \times v_1 \rangle_{\partial\kappa}. \label{eq:midway_3}
\end{equation}
Next we set $v_3 = \nabla(\nabla \times v_1)$; substituting \cref{eq:midway_2} into \cref{eq:midway_3} 
and summing over all elements in the mesh $\mathcal{T}^h$ gives
\begin{align}
\sum_{\kappa \in \mathcal{T}^h} (\nabla \times \vec{f}, v_1)_\kappa &= \sum_{\kappa \in \mathcal{T}^h} (\nabla(\nabla \times \phi), G^\top \nabla (\nabla \times v_1))_\kappa - \sum_{\kappa \in \mathcal{T}^h} (P \mathbb{I}, \nabla (\nabla \times v_1))_\kappa \nonumber \\
&\phantom{=}+ \sum_{\kappa \in \mathcal{T}^h} \langle \widehat{\phi} - \phi, \vec{n} \times (\nabla \cdot (G^\top \nabla (\nabla \times v_1))) \rangle_{\partial\kappa} \nonumber \\
&\phantom{=}+ \sum_{\kappa \in \mathcal{T}^h} \langle (\widehat{F}_3 - \nabla \times \phi) \otimes \vec{n}, G^\top \nabla (\nabla \times v_1) \rangle_{\partial\kappa} \nonumber \\
&\phantom{=}- \sum_{\kappa \in \mathcal{T}^h} \langle \widehat{F}_2, (\nabla \times v_1) \otimes \vec{n} \rangle_{\partial\kappa} \nonumber \\
&\phantom{=}- \sum_{\kappa \in \mathcal{T}^h} \langle \widehat{F}_1, \vec{n} \times v_1 \rangle_{\partial\kappa}. \label{eq:primal_form_pre}
\end{align}
Finally, given
\begin{equation}
(P \mathbb{I}, \nabla (\nabla \times v_1))_\kappa = (P, \nabla \cdot (\nabla \times v_1))_\kappa = 0
\end{equation}
and replacing the flux approximations with their corresponding Neumann
boundary data on $\partial\Omega_{N}$ and $\partial\Omega_{\phi,N}$,
as well as the analytical solution with the (DG) finite element approximation
$\phi_h \in V^{h,p}_\text{DG} = \{ v \in L_2(\Omega) : \left. v \right|_{\kappa} \in \mathcal{P}_p(\kappa) \:\: \forall \kappa \in \mathcal{T}^h\}$, $p\geq 2$, and selecting $v_1 = \psi$ gives rise to the flux
formulation: find $\phi_h \in V^{h,p}_\text{DG}$ such that
\begin{align}
\sum_{\kappa \in \mathcal{T}^h} (\nabla \times \vec{f}, \psi)_\kappa &= \sum_{\kappa \in \mathcal{T}^h} (\nabla(\nabla \times \phi_h), G^\top \nabla (\nabla \times \psi))_\kappa \nonumber \\
&\phantom{=}- \sum_{\kappa \in \mathcal{T}^h} \langle \vec{f} \times \vec{n}, \psi \rangle_{\partial\kappa\cap\partial\Omega_{\phi,N}}
- \sum_{\kappa \in \mathcal{T}^h} \langle \vec{g}_{N}, \nabla \times \psi \rangle_{\partial\kappa\cap\partial\Omega_{N}} \nonumber \\
&\phantom{=}+ \sum_{\kappa \in \mathcal{T}^h} \langle \widehat{\phi}_h - \phi_h, \vec{n} \times (\nabla \cdot (G^\top \nabla (\nabla \times \psi))) \rangle_{\partial\kappa\setminus\partial\Omega_{\phi,N}} \nonumber \\
&\phantom{=}+ \sum_{\kappa \in \mathcal{T}^h} \langle (\widehat{F}_3 - \nabla \times \phi_h) \otimes \vec{n}, G^\top \nabla (\nabla \times \psi) \rangle_{\partial\kappa\setminus\partial\Omega_{N}} \nonumber \\
&\phantom{=}- \sum_{\kappa \in \mathcal{T}^h} \langle \widehat{F}_2, (\nabla \times \psi) \otimes \vec{n} \rangle_{\partial\kappa\setminus\partial\Omega_{N}} \nonumber \\
&\phantom{=}- \sum_{\kappa \in \mathcal{T}^h} \langle \widehat{F}_1, \vec{n} \times \psi \rangle_{\partial\kappa\setminus\partial\Omega_{\phi,N}}
\quad \forall \psi \in V^{h,p}_\text{DG}. \label{eq:primal_form}
\end{align}

\subsection{Primal formulation}

We define the specialized average operators
\begin{equation}
\left.\avga{\cdot}\right|_F = \ell^+ (\cdot)^+ + \ell^- (\cdot)^- 
\text{ and } 
\left.\avgmp{\cdot}\right|_F = \ell^- (\cdot)^+ + \ell^+ (\cdot)^- \quad F \in \Gamma_I,
\end{equation}
where the weights are constrained by $\ell^+ + \ell^- = 1$; the relationship to the weights $w^{\pm}$ will become evident below.
For sufficiently smooth vector and scalar valued functions $\vec{v}$ and $q$, respectively, we define the following jump operators
\begin{align*}
\left.\jump{\vec{v}}\right|_F   &= \vec{n}^+ \cdot \vec{v}^+ + \vec{n}^- \cdot \vec{v}^-, & F \in \Gamma_I, \\
\left.\tjump{\vec{v}}\right|_F  &= \vec{n}^+ \times \vec{v}^+ + \vec{n}^- \times \vec{v}^-, & F \in \Gamma_I, \\
\left.\jump{q}\right|_F   &= \vec{n}^+ q^+ + \vec{n}^- q^-, & F \in \Gamma_I. \\
\end{align*}
\begin{lemma}
For sufficiently smooth $q : \Omega \mapsto \mathbb{R}$,
$\vec{v}, \vec{z} : \Omega \mapsto
\mathbb{R}^2$ and $\sigma : \Omega \mapsto \mathbb{R}^{2\times2}$, the following identities hold
\begin{align}
\sum_{\kappa \in \mathcal{T}^h} \langle q, \vec{v} \cdot \vec{n} \rangle_{\partial\kappa}
&= \langle \jump{q}, \avgmp{\vec{v}} \rangle_{\Gamma_I}
+ \langle \avga{q}, \jump{\vec{v}} \rangle_{\Gamma_I}
+ \langle q \vec{n}, \vec{v} \rangle_{\partial\Omega},
\label{eq:hard_to_prove1} \\
\sum_{\kappa \in \mathcal{T}^h} \langle \sigma, \vec{v} \otimes \vec{n} \rangle_{\partial\kappa}
&= \langle \ojump{\vec{v}}, \avgmp{\sigma} \rangle_{\Gamma_I}
+ \langle \avga{\vec{v}}, \jump{\sigma} \rangle_{\Gamma_I}
+ \langle \vec{v} \otimes \vec{n}, \sigma \rangle_{\partial\Omega},
\label{eq:hard_to_prove2}  \\
\sum_{\kappa \in \mathcal{T}^h} \langle \vec{z}, \vec{n} \times \vec{v} \rangle_{\partial\kappa}
&= \langle \tjump{\vec{v}}, \avgmp{\vec{z}} \rangle_{\Gamma_I}
- \langle \avga{\vec{v}}, \tjump{\vec{z}} \rangle_{\Gamma_I}
+ \langle \vec{n} \times \vec{v}, \vec{z} \rangle_{\partial\Omega}.
\label{eq:hard_to_prove3} 
\end{align}
\end{lemma}
\begin{proof}
Consider \cref{eq:hard_to_prove1}. In the absence of interior mesh facets, on the boundary it is clear that \cref{eq:hard_to_prove1} trivially holds, i.e.,
\begin{equation}
\sum_{\kappa \in \mathcal{T}^h} \langle q, \vec{v} \cdot \vec{n} \rangle_{\partial\kappa \cap \partial\Omega}
= \langle q \vec{n}, \vec{v} \rangle_{\partial\Omega}. \nonumber
\end{equation}
Let us now consider the interior mesh facets; to this end the left-hand side of \cref{eq:hard_to_prove1} may be rewritten as follows
\begin{align*}
\sum_{\kappa \in \mathcal{T}^h} &\langle q, \vec{v} \cdot \vec{n} \rangle_{\partial\kappa \setminus \partial\Omega} = \langle q^+, \vec{v}^+ \cdot \vec{n}^+ \rangle_{\Gamma_I} + \langle q^-, \vec{v}^- \cdot \vec{n}^- \rangle_{\Gamma_I}. \nonumber \\
\end{align*}
Using the fact that $1 = \ell^+ + \ell^-$, we get
\begin{align*}
\sum_{\kappa \in \mathcal{T}^h} \langle q, \vec{v} \cdot \vec{n} \rangle_{\partial\kappa \setminus \partial\Omega}
&= \langle \ell^+ q^+, \vec{v}^+ \cdot \vec{n}^+ \rangle_{\Gamma_I} + \langle \ell^- q^+, \vec{v}^+ \cdot \vec{n}^+ \rangle_{\Gamma_I} 
+ \langle \ell^+ q^-, \vec{v}^- \cdot \vec{n}^- \rangle_{\Gamma_I} \\
&\phantom{=}  + \langle \ell^- q^-, \vec{v}^- \cdot \vec{n}^- \rangle_{\Gamma_I}. \nonumber \\
\end{align*}
Noting that on an interior mesh facet we have
 $q^+ \vec{v}^- \cdot \vec{n}^+ + q^+ \vec{v}^- \cdot \vec{n}^- = 0$ and $q^- \vec{v}^+ \cdot \vec{n}^+ + q^- \vec{v}^+ \cdot \vec{n}^- = 0$, then
\begin{align*}
\sum_{\kappa \in \mathcal{T}^h} &\langle q, \vec{v} \cdot \vec{n} \rangle_{\partial\kappa \setminus \partial\Omega} \\
&= \langle \ell^+ q^+, \vec{v}^+ \cdot \vec{n}^+ \rangle_{\Gamma_I} + \langle \ell^- q^+, \vec{v}^+ \cdot \vec{n}^+ \rangle_{\Gamma_I} 
+ \langle \ell^+ q^-, \vec{v}^- \cdot \vec{n}^- \rangle_{\Gamma_I} + \langle \ell^- q^-, \vec{v}^- \cdot \vec{n}^- \rangle_{\Gamma_I} \nonumber \\
&\phantom{=} + \langle \ell^+ q^+, \vec{v}^- \cdot \vec{n}^+ \rangle_{\Gamma_I} + \langle \ell^- q^+, \vec{v}^- \cdot \vec{n}^+ \rangle_{\Gamma_I} + \langle \ell^+ q^+, \vec{v}^- \cdot \vec{n}^- \rangle_{\Gamma_I} + \langle \ell^- q^+, \vec{v}^- \cdot \vec{n}^- \rangle_{\Gamma_I} \nonumber \\
&\phantom{=} + \langle \ell^+ q^-, \vec{v}^+ \cdot \vec{n}^+ \rangle_{\Gamma_I} + \langle \ell^- q^-, \vec{v}^+ \cdot \vec{n}^+ \rangle_{\Gamma_I} + \langle \ell^+ q^-, \vec{v}^+ \cdot \vec{n}^- \rangle_{\Gamma_I} + \langle \ell^- q^-, \vec{v}^+ \cdot \vec{n}^- \rangle_{\Gamma_I}. \nonumber \\
\end{align*}
Collecting coefficients of $\vec{v}^+ \cdot \vec{n}^+$ and $\vec{v}^- \cdot \vec{n}^-$ gives
\begin{align*}
\sum_{\kappa \in \mathcal{T}^h} &\langle q, \vec{v} \cdot \vec{n} \rangle_{\partial\kappa \setminus \partial\Omega}
= \langle \avga{q}, \jump{\vec{v}} \rangle_{\Gamma_I} \nonumber \\
&\phantom{=}+ \langle \ell^- q^+, \vec{v}^+ \cdot \vec{n}^+ \rangle_{\Gamma_I} 
+ \langle \ell^+ q^-, \vec{v}^- \cdot \vec{n}^- \rangle_{\Gamma_I} \nonumber \\
&\phantom{=} + \langle \ell^+ q^+, \vec{v}^- \cdot \vec{n}^+ \rangle_{\Gamma_I} + \langle \ell^- q^+, \vec{v}^- \cdot \vec{n}^+ \rangle_{\Gamma_I} + \langle \ell^- q^+, \vec{v}^- \cdot \vec{n}^- \rangle_{\Gamma_I} \nonumber \\
&\phantom{=} + \langle \ell^+ q^-, \vec{v}^+ \cdot \vec{n}^+ \rangle_{\Gamma_I} + \langle \ell^+ q^-, \vec{v}^+ \cdot \vec{n}^- \rangle_{\Gamma_I} + \langle \ell^- q^-, \vec{v}^+ \cdot \vec{n}^- \rangle_{\Gamma_I}. \nonumber \\
\end{align*}
If we now collect coefficients of $q^+ \vec{n}^+$ and $q^- \vec{n}^-$, we get
\begin{align*}
\sum_{\kappa \in \mathcal{T}^h} &\langle q, \vec{v} \cdot \vec{n} \rangle_{\partial\kappa \setminus \partial\Omega}
= \langle \avga{q}, \jump{\vec{v}} \rangle_{\Gamma_I} + \langle \avgmp{\vec{v}}, \jump{q} \rangle_{\Gamma_I} \nonumber \\
&\phantom{=} + \langle \ell^- q^+, \vec{v}^- \cdot \vec{n}^+ \rangle_{\Gamma_I} + \langle \ell^- q^+, \vec{v}^- \cdot \vec{n}^- \rangle_{\Gamma_I}
 + \langle \ell^+ q^-, \vec{v}^+ \cdot \vec{n}^+ \rangle_{\Gamma_I} + \langle \ell^+ q^-, \vec{v}^+ \cdot \vec{n}^- \rangle_{\Gamma_I}, \nonumber \\
&=\langle \avga{q}, \jump{\vec{v}} \rangle_{\Gamma_I} + \langle \avgmp{\vec{v}}, \jump{q} \rangle_{\Gamma_I}. \nonumber
\end{align*}
Noting that $\langle \vec{z}, \vec{n} \times \vec{v}\rangle_{\partial\kappa\setminus\partial\Omega} = 
-\langle \vec{n} \times \vec{z}, \vec{v}\rangle_{\partial\kappa\setminus\partial\Omega}$
\cref{eq:hard_to_prove2,eq:hard_to_prove3} follow analogously.
\end{proof}
Employing the identities in \cref{eq:hard_to_prove1,eq:hard_to_prove2,eq:hard_to_prove3}, \cref{eq:primal_form} may be rewritten in the following equivalent manner
\begin{align}
(\nabla_h &\times \vec{f}, \psi)_\Omega \nonumber \\
&= (\nabla_h(\nabla_h \times \phi_h), G^\top \nabla_h (\nabla_h \times \psi))_\Omega
- \langle \vec{f} \times \vec{n}, \psi \rangle_{\partial\Omega_{\phi,N}}
- \langle \vec{g}_{N}, \nabla \times \psi \rangle_{\partial\Omega_{N}} \nonumber \\
&\phantom{=}+ \langle \avgmp{\widehat{\phi}_h - \phi_h}, \tjump{\nabla_h \cdot (G^\top \nabla_h (\nabla_h \times \psi))} \rangle_{\Gamma_I} \nonumber \\
&\phantom{=}\phantom{=}- \langle \tjump{\widehat{\phi}_h - \phi_h}, \avga{\nabla_h \cdot (G^\top \nabla_h (\nabla_h \times \psi))} \rangle_{\Gamma_I} \nonumber \\
&\phantom{=}\phantom{=}+ \langle \widehat{\phi}_h - \phi_h, \vec{n} \times \nabla_h \cdot (G^\top \nabla_h (\nabla_h \times \psi)) \rangle_{\partial\Omega_{\phi,D}} \nonumber \\
&\phantom{=}+ \langle \ojump{\widehat{F}_3 - \nabla_h \times \phi_h}, \avgmp{G^\top \nabla_h (\nabla_h \times \psi)} \rangle_{\Gamma_I} \nonumber \\
&\phantom{=}\phantom{=}+ \langle \avga{\widehat{F}_3 - \nabla_h \times \phi_h}, \jump{G^\top \nabla_h (\nabla_h \times \psi)} \rangle_{\Gamma_I} \nonumber \\
&\phantom{=}\phantom{=}+ \langle (\widehat{F}_3 - \nabla_h \times \phi_h)\otimes \vec{n}, G^\top \nabla_h (\nabla_h \times \psi) \rangle_{\partial\Omega_{D}} \nonumber \\
&\phantom{=}- \langle \avgmp{\widehat{F}_2}, \ojump{\nabla_h \times \psi} \rangle_{\Gamma_I}
- \langle \jump{\widehat{F}_2}, \avga{\nabla_h \times \psi} \rangle_{\Gamma_I}
- \langle \widehat{F}_2, (\nabla_h \times \psi)\otimes \vec{n} \rangle_{\partial\Omega_{D}} \nonumber \\
&\phantom{=}- \langle \avgmp{\widehat{F}_1}, \tjump{\psi} \rangle_{\partial\Omega_{\phi,D}}
+ \langle \tjump{\widehat{F}_1}, \avga{\psi} \rangle_{\Gamma_I}
- \langle \widehat{F}_1, \vec{n} \times \psi \rangle_{\partial\Omega_{\phi,D}}. \label{eq:primal_facet}
\end{align}

\subsection{\Gls{SIPG} formulation}

The choice of the flux approximations profoundly affects the numerical
properties of the method \citep{unifieddg}. We choose the flux approximations
$\widehat{F_1}$, $\widehat{F_2}$ and $\widehat{F_3}$ and $\widehat{\phi}_h$
according to the \gls{SIPG} formulation such that
\begin{align}
\widehat{F}_1 &=
\begin{cases}
\avga{F_1(\phi_h)} - \alpha \tjump{\phi_h} & \text{on } \Gamma_I, \\
F_{1}(\phi) - \alpha \vec{n} \times (\phi_h - \phi_D) & \text{on } \partial\Omega_{\phi,D}, \\
\end{cases} \label{eq:sipg_flux_1} \\
\widehat{F}_2 &= 
\begin{cases}
\avgmp{F_2(\nabla \times \phi_h)} - \beta \ojump{\nabla \times \phi_h} & \text{on } \Gamma_I, \\
F_{2}(\nabla \times \phi_h) - \beta (\nabla \times \phi_h - \vec{u}_D) \otimes \vec{n} & \text{on } \partial\Omega_{D},  \\
\end{cases} \label{eq:sipg_flux_2} \\
\widehat{F}_3 &=
\begin{cases}
\avga{F_3(\nabla \times \phi_h)} & \text{on } \Gamma_I, \\
F_{3}(\vec{u}_D) & \text{on } \partial\Omega_{D}, \\
\end{cases} \label{eq:sipg_flux_3} \\
\widehat{\phi}_h &=
\begin{cases}
\avgmp{\phi_h} & \text{on } \Gamma_I, \\
\phi_{D} & \text{on } \partial\Omega_{\phi,D}, \\
\end{cases} \label{eq:sipg_flux_4} 
\end{align}
Here, the terms $\alpha$ and $\beta$ are penalty parameters required
to ensure stability of the underlying method.
Inserting the flux functions in
\cref{eq:sipg_flux_1,eq:sipg_flux_2,eq:sipg_flux_3,eq:sipg_flux_4} into
\cref{eq:primal_facet} we arrive at the facet oriented formulation: find
$\phi_h \in V^{h,p}_\text{DG}$ such that
\begin{align}
(\nabla_h \times \vec{f}, \psi)_\Omega &= (\nabla_h(\nabla_h \times \phi_h), G^\top \nabla_h (\nabla_h \times \psi))_\Omega \nonumber \\
&\phantom{=}- \langle \vec{f} \times \vec{n}, \psi \rangle_{\partial\Omega_{\phi,N}}
- \langle \vec{g}_{N}, \nabla_h \times \psi \rangle_{\partial\Omega_{N}} \nonumber \\
&\phantom{=}+ \langle \phi_D - \phi_h, \vec{n} \times \nabla_h \cdot (G^\top \nabla_h (\nabla_h \times \psi)) \rangle_{\partial\Omega_{\phi,D}}
+ \langle \tjump{\phi_h}, \avga{\nabla_h \cdot (G^\top \nabla_h (\nabla_h \times \psi))} \rangle_{\Gamma_I} \nonumber \\
&\phantom{=} + \langle (\vec{u}_D - \nabla_h \times \phi_h)\otimes \vec{n}, G^\top \nabla_h (\nabla_h \times \psi) \rangle_{\partial\Omega_{D}}
- \langle \ojump{\nabla_h \times \phi_h}, \avgmp{G^\top \nabla_h (\nabla_h \times \psi)} \rangle_{\Gamma_I} \nonumber \\
&\phantom{=}- \langle \avgmp{2 \mu \varepsilon_h(\nabla_h \times \phi_h)} - \beta \ojump{\nabla_h \times \phi_h}, \ojump{\nabla_h \times \psi} \rangle_{\Gamma_I} \nonumber \\
&\phantom{=}- \langle 2 \mu \varepsilon_h(\nabla_h \times \phi_h) - \beta (\nabla_h \times \phi_h - \vec{u}_D) \otimes \vec{n}, (\nabla_h \times \psi) \otimes \vec{n} \rangle_{\partial\Omega_{D}}  \nonumber \\
&\phantom{=}- \langle \avga{- \nabla_h \cdot (2 \mu \varepsilon_h(\nabla_h \times \phi_h))} - \alpha \tjump{\phi_h}, \tjump{\psi} \rangle_{\Gamma_I} \nonumber \\
&\phantom{=}- \langle - \nabla_h \cdot (2 \mu \varepsilon_h(\nabla_h \times \phi_h)) - \alpha \vec{n} \times (\phi_h - \phi_D), \vec{n} \times \psi \rangle_{\partial \Omega_{\phi,D}}
\quad \forall \psi \in V^{h,p}_\text{DG}. \label{eq:facet_formulation}
\end{align}

\subsection{$C^0$-\gls{IPG} formulation}
\label{sec:appendix_c0ipg}

In order to reduce \cref{eq:facet_formulation} to the $C^0$-\gls{IPG}
formulation we restrict the space in which the \gls{FE} solution is sought to
$V^{h,p}_0$. Note that given the $C^0$ continuity of a function $\psi_h
\in V^{h,p}_0$ we have $(\ojump{\psi_h})_{ij} = 0$, $\tjump{\psi_h} =
\vec{0}$ and $\left.\psi_h\right|_{\partial\Omega_D} = 0$.
Furthermore we notice that for isotropic viscosity we have $2 \mu
\varepsilon(\vec{u}) = G \nabla \vec{u} = G^\top \nabla \vec{u}$, where
\begin{equation}
G = \mu
\begin{pmatrix}
\begin{pmatrix}
2 & 0 \\
0 & 0
\end{pmatrix} &
\begin{pmatrix}
0 & 1 \\
1 & 0
\end{pmatrix} \\
\begin{pmatrix}
0 & 1 \\
1 & 0
\end{pmatrix} &
\begin{pmatrix}
0 & 0 \\
0 & 2
\end{pmatrix}
\end{pmatrix},
\end{equation}
and integrating by parts the left side we have
\begin{equation}
(\nabla_h \times \vec{f}, \psi) = (\vec{f}, \nabla_h \times \psi) 
- \langle \vec{f} \times \vec{n}, \psi \rangle_{\partial\Omega_{\phi,N}}.
\end{equation}
%
Eliminating these terms from the facet oriented formulation in
\cref{eq:facet_formulation}, setting $w^\pm = \ell^\mp$ (i.e., $\avg{\cdot} =
\avgmp{\cdot}$), substituting for $\vec{u}_h = \nabla_h \times
\phi_h$ and $\vec{v}_h = \nabla_h \times \psi_h$ and rearranging for bilinear
and linear components we arrive at the formulation in
\cref{eq:bilinear_fem,eq:bilinear_B,eq:bilinear_l}.

\printnoidxglossary[type=\acronymtype,nonumberlist,nogroupskip,nopostdot]

\raggedright 
\bibliography{references}

\end{document}